\input harvmac
\input epsf
%%%%%%%%%%%%%  DEFINITIONS %%%%%%%%%%%%%%%%%%%%%%%%%%%%%%%%%%%
\def\figin{\epsfcheck\figin}\def\figins{\epsfcheck\figins}
\def\epsfcheck{\ifx\epsfbox\UnDeFiNeD
\message{(NO epsf.tex, FIGURES WILL BE IGNORED)}
\gdef\figin##1{\vskip2in}\gdef\figins##1{\hskip.5in}% blank space instead
\else\message{(FIGURES WILL BE INCLUDED)}%
\gdef\figin##1{##1}\gdef\figins##1{##1}\fi}
\def\DefWarn#1{}
\def\figinsert{\goodbreak\topinsert}
\def\ifig#1#2#3#4{\DefWarn#1\xdef#1{fig.~\the\figno}
\writedef{#1\leftbracket fig.\noexpand~\the\figno}%
\figinsert\figin{\centerline{\epsfxsize=#3mm \epsfbox{#2}}}
\bigskip\medskip\centerline{\vbox{\baselineskip12pt
\advance\hsize by -1truein\noindent\footnotefont{\sl Fig.~\the\figno:}\sl\
#4}}
\bigskip\endinsert\noindent\global\advance\figno by1}

\lref\Segalone{G.~Segal, ``Two-dimensional conformal field theories and 
modular functors,'' IXth International Congress on Mathematical Physics (Swansea, 1988), 
pp. 22 - 37.} 
\lref\Segaltwo{G.~Segal, ``The definition
of conformal field theory,''
in Topology, geometry and quantum field
theory, London Math. Soc. Lecture Note Ser.,
pp. 421 - 577, First circulated in 1988. }

\lref\NV{A.~Neitzke and C.~Vafa,
``${\cal N} = 2$ strings and the twistorial Calabi-Yau,''
{\tt hep-th/0402128}.
%%CITATION = HEP-TH 0402128;%%
}
\lref\witr{E.~Witten,
``Perturbative gauge theory as a string theory in twistor space,''
{\tt hep-th/0312171}.
%%CITATION = HEP-TH 0312171;%%
}
\lref\orv{A. Okounkov, N. Reshetikhin, and C. Vafa, ``Quantum Calabi-Yau and
classical crystals'', {\tt hep-th/0309208}.}
\lref\inov{
A.~Iqbal, N.~Nekrasov, A.~Okounkov, and C.~Vafa,
``Quantum foam and topological strings,''
{\tt hep-th/0312022}.
%%CITATION = HEP-TH 0111068;%%
}
\lref\minet{M. Aganagic, R. Dijkgraaf, A. Klemm, M. Marino, and
C. Vafa, ``Topological strings and integrable hierarchies,''
{\tt hep-th/0312085}.}
\lref\don{S.K. Donaldson and R. Friedman, ``Connected sums of self dual
manifolds and deformations of singular spaces,'' Nonlinearity {\bf 2}, 197
(1989).}
\lref\dt{S. Donaldson and R. Thomas, ``Gauge theory in higher
dimensions,'' in {\it The geometric universe: science, geometry
and the work of Roger Penrose}, S. Huggett et. al. eds., Oxford
Univ. Press, 1998.}
\lref\pandet{D. Maulik, N. Nekrasov, A. Okounkov, and R. Pandharipande,
``Gromov-Witten theory and Donaldson-Thomas theory,''
{\tt math.AG/0312059}.}
\lref\gova{
R.~Gopakumar and C.~Vafa,
``M-theory and topological strings. I,''
{\tt hep-th/9809187};
%%CITATION = HEP-TH 9809187;%%
``M-theory and topological strings. II,''
{\tt hep-th/9812127}.
%%CITATION = HEP-TH 9812127;%%
}
\lref\asl{A.~S.~Losev, ``Perspectives of string theory,''
talk at the ``String theory at Greater Paris'' seminar, 2001.}
\lref\cv{C. Vafa, work in progress.}
\lref\ahs{M.F. Atiyah, N.J. Hitchin, and I.M. Singer, ``Self-duality
in four dimensional Riemannian geometry,'' Proc. Roy. Soc. London Ser. A
{\bf 362}, 425 (1978).}
\lref\taub{C.H. Taubes, ``The Existence of anti-self-dual conformal
structures,'' J. Diff. Geom. {\bf 36}, 163 (1992).}
\lref\bcov{M. Bershadsky, S. Cecotti, H. Ooguri, and C. Vafa, ``Kodaira-Spencer
theory of gravity and exact results for quantum string amplitudes,''
Commun. Math. Phys. {\bf 165} (1994) 311-428, {\tt hep-th/9309140}.}
\lref\naret{I. Antoniadis, E. Gava, K. S. Narain, and
 T. R. Taylor, ``topological
amplitudes in string theory,''
Nucl. Phys. {\bf B413} (1994) 162-184, {\tt hep-th/9307158}.}
\lref\dew{
G.~Lopes Cardoso, B.~de Wit, and T.~Mohaupt,
``Deviations from the area law for supersymmetric black holes,''
Fortsch.\ Phys.\  {\bf 48}, 49 (2000),
{\tt hep-th/9904005}.
%%CITATION = HEP-TH 9904005;%%
}
\lref\osv{H. Ooguri, A. Strominger, and C. Vafa, work to appear.}
\lref\vaug{C. Vafa, ``Superstrings and topological strings at
large $N$,''
J.\ Math.\ Phys.\  {\bf 42}, 2798 (2001),
{\tt hep-th/0008142}.
%%CITATION = HEP-TH 0008142;%%
}
\lref\kachet{S.~Kachru, M.~B.~Schulz, P.~K.~Tripathy, and S.~P.~Trivedi,
``New supersymmetric string compactifications,''
JHEP {\bf 0303}, 061 (2003),
{\tt hep-th/0211182}.
%%CITATION = HEP-TH 0211182;%%
}
\lref\louisone{
S.~Gurrieri, J.~Louis, A.~Micu, and D.~Waldram,
``Mirror symmetry in generalized Calabi-Yau compactifications,''
Nucl.\ Phys.\ B {\bf 654}, 61 (2003),
{\tt hep-th/0211102}.
%%CITATION = HEP-TH 0211102;%%
}
\lref\louistwo{
T.W. Grimm and J. Louis, ``The effective action of
${\cal N}=1$ Calabi-Yau orientifolds,''
{\tt hep-th/0403067}.
}
\lref\gvw{
S.~Gukov, C.~Vafa, and E.~Witten,
``CFT's from Calabi-Yau four-folds,''
Nucl.\ Phys.\ B {\bf 584}, 69 (2000)
[Erratum-ibid.\ B {\bf 608}, 477 (2001)],
{\tt hep-th/9906070}.
%%CITATION = HEP-TH 9906070;%%
}
\lref\tv{
T.~R.~Taylor and C.~Vafa,
``RR flux on Calabi-Yau and partial supersymmetry breaking,''
Phys.\ Lett.\ B {\bf 474}, 130 (2000),
{\tt hep-th/9912152}.
%%CITATION = HEP-TH 9912152;%%
}
%\StromingerNS
\lref\StromingerNS{
A.~Strominger,
%``Vacuum Topology And Incoherence In Quantum Gravity,''
Phys.\ Rev.\ Lett.\  {\bf 52}, 1733 (1984).
%%CITATION = PRLTA,52,1733;%%
}
\lref\withcs{E. Witten, ``Chern-Simons gauge theory as a
string theory,'' Prog. Math. {\bf 133}, 637 (1995),
{\tt hep-th/9207094}.}

\lref\lmn{A.~Losev, A.~Marshakov, and N.~Nekrasov,
``Small instantons, little strings and free
fermions,'' {\tt hep-th/0302191},
in: {\it From fields to strings: circumnavigating theoretical physics,}
Ian Kogan Memorial Volume,
M. Shifman, A. Vainshtein, and J. Wheater eds., World
Scientific, Singapore.}

%%%% Nikita's defs %%%%%%%%
%%%%%
\def\boxit#1{\vbox{\hrule\hbox{\vrule\kern8pt
\vbox{\hbox{\kern8pt}\hbox{\vbox{#1}}\hbox{\kern8pt}}
\kern8pt\vrule}\hrule}}
\def\mathboxit#1{\vbox{\hrule\hbox{\vrule\kern8pt\vbox{\kern8pt
\hbox{$\displaystyle #1$}\kern8pt}\kern8pt\vrule}\hrule}}

%%%%%%%%
%%%%%%%%%%%% Greek %%%%%%%%%%%%

\def\u{{\Upsilon}}

%%%%%
%%%%%%%Russian fonts%%%%%%%5
\chardef\tempcat=\the\catcode`\@ \catcode`\@=11
\def\cyracc{\def\u##1{\if \i##1\accent"24 i%
    \else \accent"24 ##1\fi }}
\newfam\cyrfam

%%%%%%%%%%%%%%%% Calligraphic letters  %%%%%%%%%%%%%

%%%%%%%%%%%% Derivatives  %%%%%%%%%%%

%%%%%%%%%%% letters with bar %%%%%%%%

%%%%%%%%%% Math symbols %%%%%%%%%%%%%

%%%%%%%%%%%% end of Nikita's defs %%%%%%%%%%%%%%%
\noblackbox
\newcount\figno
 \figno=1
 \def\fig#1#2#3{
 \par\begingroup\parindent=0pt\leftskip=1cm\rightskip=1cm\parindent=0pt
 \baselineskip=11pt
 \global\advance\figno by 1
 \midinsert
 \epsfxsize=#3
 \centerline{\epsfbox{#2}}
 \vskip 12pt
 {\bf Fig.\ \the\figno: } #1\par
 \endinsert\endgroup\par
 }
 \def\figlabel#1{\xdef#1{\the\figno}}
 \def\encadremath#1{\vbox{\hrule\hbox{\vrule\kern8pt\vbox{\kern8pt
 \hbox{$\displaystyle #1$}\kern8pt}
 \kern8pt\vrule}\hrule}}
\Title
{\vbox{
 \baselineskip12pt
\hbox{CALT-68-2718}
\hbox{IPMU-09-0004}
}}
{\vbox{
 \centerline{Geometry As Seen By String Theory}
 }}
\centerline{Hirosi Ooguri}
\bigskip

\centerline{\it California Institute of Technology, Pasadena, CA 91125, USA}
\centerline{\it and}
\centerline{\it Institute for the Physics and Mathematics of the Universe,}
\centerline{\it University of Tokyo, Kashiwa, Chiba 277-8586, Japan}

\bigskip
\bigskip
\bigskip
\bigskip
\bigskip
\bigskip
\bigskip

\centerline{\sl The Fourth Takagi Lectures of
the Mathematical Society of Japan,}
\centerline{\sl  delivered on 21 June 2008 
at Department of Mathematics, Kyoto University}

\smallskip
\bigskip
\bigskip

%\vskip .1in

%\smallskip
\Date{}

\newsec{Introduction}

It is a great privilege to deliver this set of lectures in honor of 
Professor Teiji Takagi (1875 - 1960), the founding father of 
modern mathematical research 
in Japan. Professor Takagi was an alumnus of my high school, and our mathematics 
teacher liked to tell us about the local hero who realized Kronecker's 
{\it Jugendtraum} in the case of imaginary quadratic fields
by establishing Class Field Theory
\lref\kronecker{T.~Takagi, ``\"Uber eine Theorie
des relativ Abel'schen Zahlk\"orpers,'' J. of the College of Science,
Imperial Univ. of Tokyo {\bf 41}-0, 1 (1920)} \kronecker. As
a high school student, I also enjoyed reading his popular book on 
the history of modern mathematics \lref\takagi{T.~Takagi, ``Kinsei Suugaku 
Shidan,''  (Iwanami, 1995).} \takagi , where the moments of creation of new 
mathematics are vividly described. I was particularly fascinated by 
the story on elliptic functions of Gauss, Abel, and Jacobi, which turned
out to be relevant in my study of conformal field theories 
7 years later. We cannot overestimate the influence of
his legendary textbook on {\it Calculus} 
\lref\kaiseki{T.~Takagi, ``Kaiseki Gairon,''
(Iwanami, 1983).} \kaiseki\ over generations of engineers and 
scientists as well as mathematicians in Japan for the last three 
quarters of a century since it was first published. 
On my bookshelves, it stands next to 
{\it Feynman Lectures on Physics} and
{\it Landau-Lifshitz Course of Theoretical Physics}, and 
I still consult it from time to time. 

The main subject of this set of lectures is the topological string theory. 
The topological string theory was introduced by E.~Witten about 20 years ago,
and it has been developed by collaborations of physicists and 
mathematicians. Its mathematical structure is very rich, and it has 
lead to discoveries of new connections between different areas of 
mathematics, ranging from algebraic geometry, symplectic geometry 
and topology, to combinatorics, probability and representation 
theory. The topological string theory also has many important applications 
to problems in physics. Though the theory was initially thought of as 
a simple toy model of string theory, it has turned out to be useful 
in computing a certain class of 
scattering amplitudes of physical string theory. In the past 10 years, 
the relation between topological string and physical string 
has been applied to variety of problems, and it has advanced 
our understanding of string compactifications, provided a powerful
computational tool to study strongly coupled dynamics of gauge theories, 
has shed light on mysteries of quantum gravity such as quantum 
states of black holes, and pointed out a promising direction to 
prove the AdS/CFT correspondence. Moreover, the
topological string theory has given us insights into how
our concept of space and time should be modified in order 
to formulate fundamental laws of nature.

Although these lectures are meant to be for mathematicians,
I felt it would be appropriate to spend the first couple of minutes
in this course explaining physicists' motivation to
study string theory. 
In the past few hundred years, physicists have searched for 
fundamental laws of nature by exploring phenomena at shorter 
and shorter distances. Although the idea that everything on the Earth is made
of atoms goes back to Ancient Greek, the modern atomic theory
began with the publication of ``New System of Chemical Philosophy''
by J. Dalton in 1808. In the middle of the 19th century, the
size of atoms is correctly estimated to be about $10^{-10}$ meters. 
By the end of the 19th century, due to the discovery of the
electron and study of radioactivity, scientists began to think that
atoms are not fundamental and that they have internal structure. 
In 1904, H. Nagaoka proposed the model in which there is
a positively charged nucleus at the center
with electrons orbiting around them. The existence of atomic
nuclei was confirmed by the Geiger-Marsden experiment and
the theory of E. Rutherford. The radius of the atomic nucleus is
about $10^{-15}$ - $10^{-14}$ meters. In the 1930th, thanks to the discovery of 
the neutron by J. Chadwick, the splitting of the atomic nucleus
by J. Cockroft and E. Walton, and the meson theory of H. Yukawa, it became
clear that the atomic nucleus is made of protons and neutrons
bound together by the $\pi$ meson. The radius of the proton is
about $10^{-15}$ meter. The progress of elementary particle physics
in the past 50 years has culminated in the ``Standard Model of Particle
Physics,'' which describes all known particle physics phenomena down
to $10^{-18}$ meters. The Large Hadron Collider, which just began its
operation at CERN in Switzerland, will probe distance as short as
$10^{-19}$ meters.   

It is natural to ask whether this progression continues 
indefinitely. Surprisingly, there are reasons to think that 
the hierarchical structure of nature will terminate at the Planck 
length at $10^{-35}$ meters. 
Let us perform a thought-experiment to explain 
why this might be the case. Physicists build particle colliders 
to probe short distances. The more energy we use to collide 
particles, the shorter distances we can explore. This has been 
the case so far. One may then ask: can we build a collider with 
energy so high that it can probe distances shorter than the 
Planck length? The answer is no. When we collide particles with 
such high energy, a black hole will form and its event horizon 
will conceal the entire interaction area. Stated in another way, 
the measurement at this energy would perturb the geometry so much 
that the fabric of space and time would be torn apart. This would 
prevent physicists from ever seeing what is happening at distances 
shorter than the Planck length. This is a new kind of uncertainty 
principle. The Planck length is truly fundamental since it is the 
distance where the hierarchical structure of nature will terminate. 

Space and time do not exist beyond the Planck scale, and they 
should emerge from a more fundamental structure. Superstring 
theory is a leading candidate for a mathematical framework to 
describe physical phenomena at this scale
since it contains all the ingredients 
necessary to unify quantum mechanics and general relativity.

\newsec{From Points to Strings}  

The axiomatic method for geometry invented by Euclid of Alexandria
is based on mathematical points with no size and structure. 
The {\it Elements} defines a point as ``that which has no part.'' 
2300 years after Euclid, string theory is offering the first real 
alternative to this approach by introducing finite-sized objects 
as basic building blocks.

Consider a Riemannian manifold $M$ and try to probe it using a
point-like particle. A typical example of ``observables'' 
is Green's function $G(x,y)$ obeying
\eqn\green{
 (-\Delta_x + m^2) G(x,y) = \delta(x,y),}
where $x,y \in M$, $\Delta_x$ is the Laplace-Beltrami operator on $x$, 
$m$ is a constant, which physicists regard as a mass of the
particle, and $\delta(x,y)$ is the delta-function for $x=y$.
Green's function can be expressed as a path integral,
$i.e.$, a sum over all possible paths in $M$
from $x$ to $y$ with an appropriate weight.

In string theory, an analogue of Green's function has 
richer structure. An obvious analogue would be
a sum over all possible spherical surface connecting $x$ to $y$ 
as in Figure 1(a). But there is no reason to stop at two points.
We can choose $n$ point, $x_1, x_2, \cdots, x_n$, and 
sum over all spherical surfaces connecting them as in Figure 1(b). To define
a similar object in a point particle theory, one would
need to introduce
``interactions.'' In string theory, interactions are already built
in without additional assumptions. More generally, we can consider
a sum over genus-$g$ surfaces connecting the $n$ points to define
an amplitude $F_g(x_1, \cdots, x_n)$. We can go even further; since 
we are considering string theory, we should be able to consider
$n$ configurations of strings in $M$ and a sum over genus-$g$ surfaces
with $n$ boundaries connecting them as in Figure 1(c). 
This can be made a little more precise as follows.

%%%%%%%%%%%%%%%%%%%%%%%%%%%%%%%%%%%%%%%%% 
\bigskip 
\centerline{\epsfxsize 6truein\epsfbox{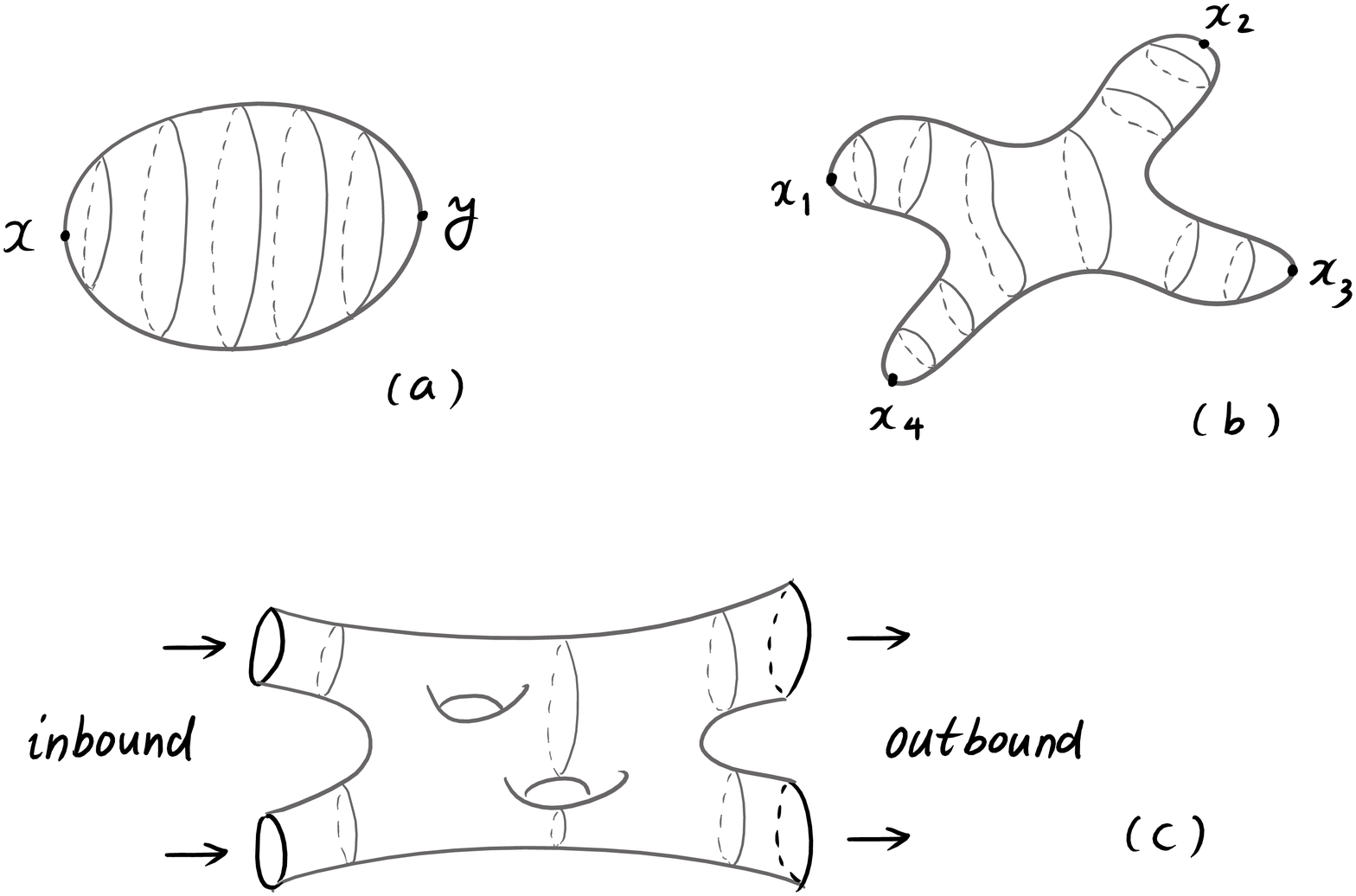}} 
%\leftskip 2pc 
%\rightskip 2pc 
\noindent{\ninepoint\sl \baselineskip=2pt {\bf 
Figure 1} {{An analogue of Green's function in string theory
is given by summing over spheres with two points fixed (a). 
This can be generalized to $n$-point functions (b) and 
to higher genus amplitudes 
with incoming and outgoing string states (c).}}} 
\bigskip 
%%%%%%%%%%%%%%%%%%%%%%%%%%%%%%%%%%%%%%%%%%%%%%%%%%%%%%% 

\subsec{String Amplitudes at Genus $g$}

Let us define a conformal field theory in two dimensions. 
We start with a Hilbert space ${\cal H}$, which realizes the Virasoro algebra
\eqn\virasoro{
[ L_n , L_m ] = (n-m) L_{n+m} + {c \over 12} (n^3 - n) \delta_{n+m, 0},
~~~ n, m = 0, \pm 1, \pm 2, \cdots.}
Physicists regard it as a ``space of states'' of the conformal field 
theory. Here $c$ is the central charge, which commutes 
with all other generators and takes a fixed value
on ${\cal H}$. The Hilbert space ${\cal H}$ is decomposed
into a sum of products of irreducible unitary representations $V_h$
of the Virasoro algebra with the highest weight $h$ as
\eqn\irreps{ {\cal H} = \oplus_{h, \bar h} N_{h,\bar h} V_h \otimes V_{\bar h},}
where $N_{h,\bar h}$ are non-negative integers.
It is convenient to distinguish the Virasoro algebra realized on the two 
factors of $V_h \oplus V_{\bar h}$, so physicists use $\{ L_n\}$ for $V_h$
and $\{ \bar{L}_n\}$ for $V_{\bar h}$ and call them the left-movers and the 
right-movers.  

According to \refs{\Segalone,\Segaltwo}, a conformal field theory 
is defined as a functor,
\eqn\functor{ \Sigma_g(n,m) \rightarrow 
{\cal A}_g(n,m) \in {\rm Hom}({\cal H}^{\otimes n}, {\cal H}^{\otimes m}),}
where $\Sigma_g(n,m)$ is a Riemann surface of genus $g$ with $n$ 
parametrized boundaries with the inbound orientation and $m$ parametrized
boundaries with the outbound orientation. 
It is supposed to satisfy the gluing axioms spelled out in \Segaltwo . 

To define string theory, we need a particular type of conformal field theory
where one can define a nilpotent operator $Q: {\cal H} \rightarrow {\cal H}$
of degree $1$ such that 
the Virasoro generators are $Q$-trivial, 
\eqn\antighost{ L_n = \{ Q, b_n \},~\bar{L}_n = \{Q, \bar{b}_n \}, }
for some operators $b_n, \bar{b}_n: 
{\cal H} \rightarrow {\cal H}$ of degree $-1$. 
We need the central charge
$c=0$ for this to be possible. In physics literature, 
$b_n$'s are called anti-ghosts. Let ${\cal M}_g(n,m)$ be the moduli space
of $\Sigma_g(n,m)$, and $\Omega^*({\cal M}_g(n,m))$ be the space of 
differential forms on it. When \antighost\ holds for a conformal field 
field theory \functor, one can define
\eqn\stringfunctor{ \omega_g(n,m) 
\in \Omega^*({\cal M}_g(n,m)) \otimes {\rm Hom}
({\cal H}^{\otimes n}, {\cal H}^{\otimes m}),}
so that it is closed, $D \omega_g(n,m) = 0$ with respect to $D = d + Q$,
where $d$ is the de Rham operator on $\Omega^*({\cal M}_g(n,m))$. 
There is a set of gluing axioms for $\omega_g(n,m)$. For example, 
for the gluing map, 
\eqn\glueone{ {\sl glue}:  {\cal M}_{g_1}(n,k)
\oplus {\cal M}_{g_2}(k,m)\rightarrow {\cal M}_{g_1+g_2}(n,m),}
$\omega_g(n,m)$ responds as
\eqn\gluetwo{
\eqalign{ {\sl glue}^*
&\left( \omega_{g_1+g_2}(n,m)\right)
= {\sl Tr}_{{\cal H}^{\otimes  k}}\left( \omega_{g_1}(n,k)
\omega_{g_2}(k,m) \right) \cr
&~~~~~~~~~~~~~~~\in \Omega^*({\cal M}_{g_1}(n,k)
\oplus {\cal M}_{g_2}(k,m)) 
\otimes {\rm Hom}({\cal H}^{\otimes n}, {\cal H}^{\otimes m}).}}
Physicists have developed a method to construct $\omega_g(n,m)$ for 
all known perturbative string theories, including  bosonic 
string, superstring, heterotic string and topological 
string. The genus-$g$ string amplitude 
$F_g(n,m)$ is given by integrating the top
component of $\omega_g(n,m)$ over ${\cal M}_g(n,m)$,
\eqn\closedstring{ F_g(n,m) = \int_{{\cal M}_g(n,m)} \omega_g(n,m).}

\subsec{Open Strings}

\lref\wittenCS{
  E.~Witten,
  ``Chern-Simons gauge theory as a string theory,''
  Prog.\ Math.\  {\bf 133}, 637 (1995)
  [arXiv:hep-th/9207094].
  %%CITATION = PMTMA,133,637;%%
}

Mathematics of open strings was pioneered by K.~Fukaya 
and M.~Kontsevich before physicists (except for Witten \wittenCS ) 
realized geometric significance of D branes. 
Although it is difficult to do justice to recent progress in this 
area's mathematics in this short course of lectures, let us briefly 
introduce D branes and open strings as we will need them later.

%%%%%%%%%%%%%%%%%%%%%%%%%%%%%%%%%%%%%%%%% 
\bigskip 
\centerline{\epsfxsize 4truein\epsfbox{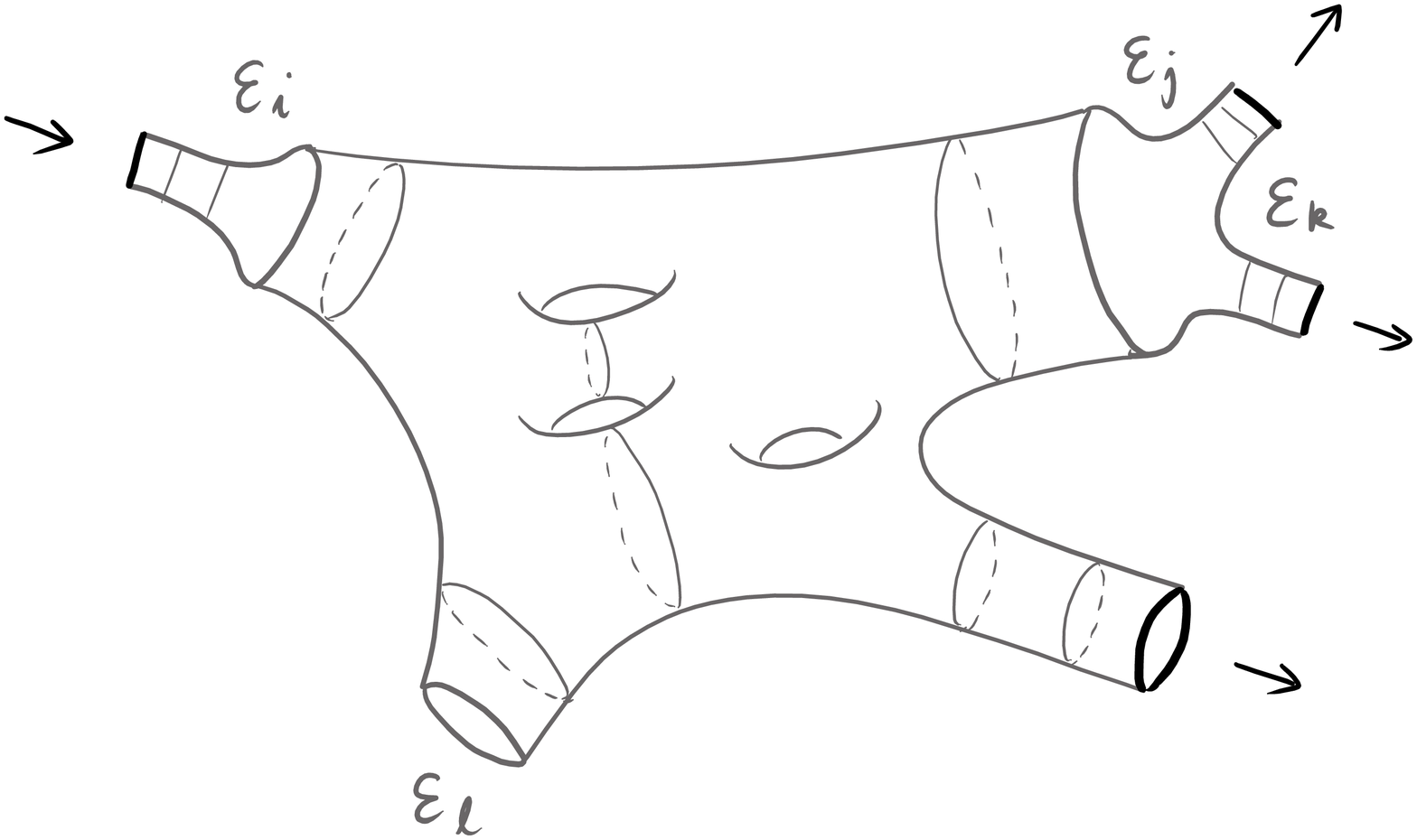}} 
%\leftskip 2pc 
%\rightskip 2pc 
\noindent{\ninepoint\sl \baselineskip=2pt {\bf 
Figure 2} {{An open string amplitude of genus 3.
It has no incoming closed string, 1 incoming open string, 
1 outgoing closed string, and 2 outgoing open strings ($n=0, n_{open}=1,
m=1, m_{open}=2$). Note that each open string state is attached to 
a segment on a boundary. The remaining boundary 
components are assigned with boundary conditions, 
${\cal E}_i, {\cal E}_j, {\cal E}_k, {\cal E}_l$.  
}}} 
\bigskip 
%%%%%%%%%%%%%%%%%%%%%%%%%%%%%%%%%%%%%%%%%%%%%%%%%%%%%%% 

An open string has two end-points, and we need to specify boundary 
conditions there. 
In physics literature, each boundary condition is called a D brane. 
Let us choose a set of D branes, $\{ {\cal E}_i: i = 1, \cdots, k \}$. 
Some of the D branes may coincide, $i.e.$, some boundary conditions
may be identical.
For a pair of D branes, ${\cal E}_i$ and ${\cal E}_j$, 
we can define an open string 
Hilbert space ${\cal H}_{open}({\cal E}_i, {\cal E}_j)$. 
Consider a Riemann surface $\Sigma_{g,b}(n,n_{open}; m, m_{open})$ with genus 
$g$ and $(b+n+m)$ boundaries with the following markings. See Figure 2. 
As in the case of closed string, there are 
$n$ parametrized boundaries with the inbound orientation and $m$
parametrized boundaries with the outbound orientation. In addition, 
there are $(n_{open}+m_{open})$ disjoint embeddings of 
the interval $[0,1]$ into the remaining
$b$ boundaries, where $n_{open}$ of them are inbound 
and $m_{open}$ are outbound.
We arrange so that boundaries make a right-angle turn 
at each end point of 
the images of $[0,1]$. This is needed for the gluing rules
to work for open strings.
The complement of these images in the $b$ boundaries has several
disjoint components, and they are either circles or intervals. 
We assign a D brane boundary condition to each of the disjoint
components. In this way, each image of 
$[0,1]$ is sandwiched by a pair of boundary components with prescribed
boundary conditions, ${\cal E}_i$ and ${\cal E}_j$ for example. 
We can then 
associate the Hilbert space ${\cal H}_{open}({\cal E}_i, {\cal E}_j)$ to
the image on $[0,1]$ on the boundary. 
The open/closed conformal field theory is a functor, 
\eqn\openfunctor{\eqalign{
\Sigma_{g, b}(n,n_{open};m,m_{open})
\rightarrow &{\cal A}_{g, b}(n,n_{open};m,m_{open})\cr
&\in {\rm Hom}\left( {\cal H}^{\otimes n} \otimes 
{\cal H}_{open}^{\otimes n_{open}} , {\cal H}^{\otimes m} \otimes 
{\cal H}_{open}^{\otimes m_{open}} \right).}}
If the Virasoro generator acting on ${\cal H}$ and ${\cal H}_{open}$ are
$Q$-trivial, one can construct 
\eqn\omenmeasure{ 
\eqalign{\omega_{g, b}(n,n_{open};m,m_{open})&
\in \Omega^*\left(\Sigma_{g,b}(n,n_{open}; m, m_{open})\right)\cr
&~~~~~\otimes {\rm Hom}\left( {\cal H}^{\otimes n} \otimes 
{\cal H}_{open}^{\otimes n_{open}} , {\cal H}^{\otimes m} \otimes 
{\cal H}_{open}^{\otimes m_{open}} \right),}}
which is closed with respect to $D = d + Q$, and 
it can be used to define a string amplitude 
$F_{g, b}(n,n_{open};m,m_{open})$
by integrating its top component over the moduli space.

Open string theory and closed string theory (the one which contains
closed strings only) seem very different. For one thing, open string 
theory can be formulated using string field theory
\lref\wittenSFT{
  E.~Witten,
  ``Noncommutative geometry and string field theory,''
  Nucl.\ Phys.\  B {\bf 268}, 253 (1986).
  %%CITATION = NUPHA,B268,253;%%
} \wittenSFT. Mathematically, it means that the moduli space
for $\Sigma_{g,b\neq 0}$ has a nice triangulation that corresponds to
a sum of Feynman diagrams. String field theories for open 
{\it topological} strings are particularly simple. In section 5,
we will discuss 
the Chern-Simons gauge theory in three dimensions \lref\wittenCS{
  E.~Witten,
  ``Chern-Simons gauge theory as a string theory,''
  Prog.\ Math.\  {\bf 133}, 637 (1995)
  [arXiv:hep-th/9207094].
  %%CITATION = PMTMA,133,637;%%
} \wittenCS\ and random matrix models 
\lref\matrix{
  R.~Dijkgraaf and C.~Vafa,
``Matrix models, topological strings, and 
supersymmetric gauge theories,''
  Nucl.\ Phys.\  B {\bf 644}, 3 (2002)
  [arXiv:hep-th/0206255];
  %%CITATION = NUPHA,B644,3;%%
  ``On geometry and matrix models,''
  Nucl.\ Phys.\  B {\bf 644}, 21 (2002)
  [arXiv:hep-th/0207106].
  %%CITATION = NUPHA,B644,21;%%
} \matrix, as examples of string field theories. 
In contrast, a string field theory to compute closed string 
amplitudes $F_g$ \closedstring\ for all $g$ is not known, except for 
the topological B-model \lref\BCOV{
  M.~Bershadsky, S.~Cecotti, H.~Ooguri and C.~Vafa,
  ``Kodaira-Spencer theory of gravity and exact results for quantum string
  amplitudes,''
  Commun.\ Math.\ Phys.\  {\bf 165}, 311 (1994)
  [arXiv:hep-th/9309140].
  %%CITATION = CMPHA,165,311;%%
} \BCOV . The topological B-model is special since it is manifest
that  the string amplitudes
$F_g$ receive contributions only from boundaries of the moduli space
${\cal M}_g$. 

\newsec{Topological String Theory}

To define the topological string theory, 
we consider a conformal field theory with
the ${\cal N}=2$ superconformal symmetry
generated by $\{L_n, G^+_n, G^-_n, J_n\}$ and the central charge $\hat{c}$.
These generators obey the commutation relations,\foot{
A CFT connoisseur may notice that the commutation relations listed here
are slightly difference from the standard ones. In fact, they are related to
by rearranging the generators, the process called {\it topological twist}.}
\eqn\sca{\eqalign{
&[L_n, L_m] = (n-m)L_{n+m}, 
~\{G_n^+, G_m^-\} = 2L_{n+m},
\cr
& [L_n, G_m^+] = -m G_{n+m}^+,~
[L_n, G_m^-] = (n-m)G_{n+m}^-,
~[J_n, G_m^\pm] = \pm G_{n+m}^\pm,\cr
&[J_n, J_m] = \hat{c}  n\delta_{n+m}, 
~ [L_n, J_m] = -m J_{n+m}.}}
We see that the Virasoro algebra with $c=0$
is a sub-algebra and that the Virasoro generators are $Q$-trivial,
$L_n = \{ Q , b_n \}$, where $Q= G_0^+ + \bar{G}_0^+$ and $b_n={1\over 2}G_n^-$.
Thus, we can use the formalism of section 2.2 to define the genus-$g$
 string amplitudes $F_g(n,m)$. 

\subsec{Two Models for One Calabi-Yau Manifold}

A prototypical example of conformal field theories
with ${\cal N}=2$ superconformal symmetry is the supersymmetric sigma-model
whose target space is a Calabi-Yau manifold $M$. 
For a given Calabi-Yau manifold, one can define two distinct sigma-models
 with different $Q$ operators, the A-model 
and the B-model \lref\wittenmirror{
  E.~Witten,
  ``Mirror manifolds and topological field theory,'' in Mirror Symmetry I, ed. S.-T. Yau, pp. 121 - 160 (American Mathematica Society, 1998) [arXiv:hep-th/9112056].
  %%CITATION = HEP-TH/9112056;%%
} \wittenmirror .

In physics, the A-model is defined in terms of a path integral over (not necessarily 
holomorphic) maps 
\eqn\map{X: \Sigma \rightarrow M,}
together with Grassmannian fields on $\Sigma$, 
\eqn\fermionsa{\eqalign{ 
 \theta(z,\bar z) &\in T_{X(z,\bar z)}^{1,0}M ,~~
\eta(z,\bar z) \in \Omega^{1,0}(\Sigma)\otimes T_{X(z,\bar z)}^{0,1}M,\cr
\bar\theta(z,\bar z)&\in T_{X(z,\bar z)}^{0,1}M , ~~
 \bar\eta(z,\bar z) \in \Omega^{0,1}(\Sigma)\otimes T_{X(z,\bar z)}^{1,0}M,}}
where $(z,\bar z)\in \Sigma$, 
$T^{1,0}_xM$ and $T^{0,1}_xM$ are the holomorphic and anti-holomorphic
components of the tangent space at $x \in M$, and $\Omega^{1,0}(\Sigma)$ and
$\Omega^{0,1}(\Sigma)$ are the spaces of $(1,0)$ and $(0,1)$ forms on $\Sigma$. 
In this model, the operator $Q$ is the Noether charge associated to the following 
transformation,
\eqn\brsta{
\eqalign{ ~&\delta X^i = \epsilon \theta^i, ~~ 
\delta \bar X^{\bar i} = \epsilon \bar \theta^{\bar i} ; ~~
          \delta \theta^i = \delta \bar \theta^{\bar i} = 0; \cr
         ~& \delta \bar \eta^i  = \epsilon \bar \partial X^i -
\Gamma^i_{jk} \theta^j \bar \eta^k, ~~
          \delta \eta^{\bar i} = \epsilon \partial\bar X^{\bar i}
- \Gamma^{\bar i}_{\bar j \bar k} \bar\theta^{\bar j} \eta^{\bar k} ,}}
where I have used holomorphic coordinates $(x^1,\cdots, x^{\hat c})$, 
$\Gamma^i_{jk}$ is the Christoffel connection associated 
to the Ricci flat K\"ahler metric, and $\epsilon$ is a Grassmann number to 
parametrize the transformation. 
It then follows that $Q$ cohomology of ${\cal H}$ in this model is 
given by the Hodge-de Rham cohomology of $M$ as,
\eqn\acohomology{ 
{\rm H}_{Q}({\cal H}) = \oplus_{p,q=0}^{\hat c} {\rm H}^{p,q}(M).}
The degrees $(p,q)$ can be identified as eigenvalues of 
$(J_0, \bar J_0)$ in the left-moving and right-moving
${\cal N}=2$ superconformal algebras.

Another consequence of the $Q$ symmetry is that a path integral for 
a $Q$ invariant amplitude localizes to a sum over fixed points of the symmetry.\foot{See
section 5 of \wittenmirror\ for discussion on the localization mechanism.}  
For the A-model, the fixed points of \brsta\ consist of holomorphic maps
from $\Sigma$ to $M$ together with $\theta = \bar \theta = 0$, so that
$\delta \bar \eta \sim  \bar \partial X = 0$. 
The space of holomorphic maps is finite dimensional 
and can be used as a basis for mathematical investigation. This leads to
the connection between the topological string theory and the Gromov-Witten 
invariants, which we will discuss later. 
An example of $Q$ invariant amplitudes is ${\cal A}_g (n,m)$
given in \functor\ evaluated for a $Q$-closed state in 
${\cal H}^{\otimes n} \otimes {\cal H}^{*\otimes m}$.
Such amplitudes can be expressed as sums over holomorphic maps. 
On the other hand, $\omega_g (n,m)$ as in \stringfunctor\ is 
invariant under $D=d+Q$ but not under $Q$. 
This leads to an interesting
subtlety in evaluating $F_g = \int_{{\cal M}_g(n,m)} \omega_g (n,m)$.
In particular, it will show up as the holomorphic anomalies
\lref\BCOVone{
  M.~Bershadsky, S.~Cecotti, H.~Ooguri and C.~Vafa,
  ``Holomorphic anomalies in topological field theories,''
  Nucl.\ Phys.\  B {\bf 405}, 279 (1993)
  [arXiv:hep-th/9302103].
  %%CITATION = NUPHA,B405,279;%%
}
\refs{\BCOVone, \BCOV}, as we will see below. 

Let us turn to the B-model. Its degrees of freedom are maps $X: \Sigma 
\rightarrow M$ and Grassmannian fields on $\Sigma$, 
\eqn\fermionsb{\eqalign{ 
\bar \theta(z,\bar z) &\in T^{0,1}_{X(z,\bar z)}M,~~\theta(z,\bar z)
\in T^{*1,0}_{X(z,\bar z)}M,
\cr
\eta(z,\bar z) &\in \left(\Omega^{1,0}(\Sigma)\oplus \Omega^{0,1}(\Sigma)\right)
\otimes T_{X(z,\bar z)}^{1,0}M.}}
The $Q$ transformation is given by
\eqn\brstb{\delta X^i = 0, ~~\delta \bar X^{\bar i} 
= \epsilon \bar\theta^{\bar i};
~~\delta \eta^i = \epsilon d X^i, ~~\delta \theta = \delta\bar \theta
=0}
The $Q$ cohomology of ${\cal H}$ in this case is the $\bar \partial$ cohomology, 
\eqn\bcohomology{ 
{\rm H}_{Q}({\cal H}) = \oplus_{p,q=0}^{\hat c}
{\rm H}_{\bar \partial}^{p}(M, \wedge^q T^{1,0} M).}
The degrees $(p,q)$ can be identified as eigenvalues of 
$(J_0, \bar J_0)$ in the left-moving and right-moving
${\cal N}=2$ superconformal algebras. 
On a Calabi-Yau manifold of complex dimensions $\hat c$, there is a unique 
holomorphic $(\hat{c},0)$-form $\Omega$ which is nowhere vanishing, and it can be
used to identify this cohomology with the Hodge-de Rham cohomology as,
\eqn\identify{{\rm H}_{\bar \partial}^{p}(M, \wedge^q T^{1,0} M)
\simeq {\rm H}^{p, \hat{c}-q}(M).}

Fixed points of the $Q$ transformation \brstb\ 
can be found at $dX = 0$, namely constant maps
$X: \Sigma \rightarrow p \in M$. Thus, $Q$ invariant amplitudes, such 
as conformal field theory amplitudes ${\cal A}_g(n, m)$ evaluated 
for $Q$ closed states in ${\cal H}^{\otimes n}
\otimes {\cal H}^{*\otimes m}$, can be expressed as a sum over constant maps, 
which is the same thing as an integral over $M$. The string amplitudes
$F_g(n,m) =\int_{{\cal M}_g(n,m)} \omega_g(n,m)$, on the other hand,
does not necessarily  
localize to integrals over $M$ since they are invariant under $D = d + Q$ and
not under $Q$. The difference of $D$ and $Q$ shows up as contributions from
boundaries of the moduli space ${\cal M}_g(n,m)$, and the part of $F_g$ 
which fails to be localized on constant maps can be expressed in terms
of Feynman diagrams for point particles. This is the origin 
of the Kodaira-Spencer description of the B-model mentioned at the end 
of section 2.2.

\subsec{Moduli Space of Topological String Theory}

According to Yau's theorem, a Ricci-flat K\"ahler metric of a Calabi-Yau manifold
$M$ is uniquely determined by complex structure and K\"ahler class.
Infinitesimal deformations of complex structure correspond to elements of
${\rm H}_{\bar \partial}^{1}(M, T^{1,0} M)$, while K\"ahler class is
parametrized by $H^{1,1}(M)$. In string theory, it is natural to complexify 
K\"ahler class. For each choice of complex structure and complexified 
K\"ahler class, there is a conformal field theory and we can use it to define
string amplitudes $F_g(n,m)$. Thus, topological string amplitudes 
should be regarded as geometric
objects over the moduli space of Calabi-Yau manifolds. 

Infinitesimal deformations of a conformal field theory are generated by 
marginal operators, which in the case of a {\it topologically twisted} 
${\cal N}=2$ conformal field theory correspond to $Q$ cohomology elements
of $(J_0, \bar J_0) = (1,1)$. In the A and B-models, they are elements 
of $H^{1,1}(M)$ and $H_{\bar \partial}^{1}(M, T^{1,0}M)$ respectively,
as we can see from \acohomology\ and \bcohomology .
Thus, we expect that string amplitudes $F_g$ depend on the K\"ahler moduli in 
the A-model and the complex moduli in the B-model. 

Due to the index theorem and the triviality of the first Chern class of $M$,
the number of zero modes of fermions $\eta,\bar \eta$ minus the 
number of zero modes of $\theta, \bar\theta$ on $\Sigma_g$
is equal to $(2g-2)\hat{c}$. When $\hat c = 3$, this coincides 
with ${\rm dim} \ {\cal M}_g$. Because of this, $\omega_g$ with $n=m=0$ is 
a top form on ${\cal M}_g$, and we can define the {\it vacuum 
amplitude} $F_g = \int_{{\cal M}_g} \omega_g$. 
We will mainly consider the case of $\hat c = 3$ 
in the following, except for $g=1$, where the index vanishes.
It turns out that this is also the most interesting case for
physical applications of string theory since $2\hat c = 6 = 10-4$, 
where $10$ is the critical dimensions of superstring theory 
and $4$ is the macroscopic dimensions of our spacetime. 

Let us discuss the moduli space of the topological string theory. 
It is easier to start with the B-model since its conformal field theory
amplitudes are expressed as integrals over $M$ and we can use classical 
geometry to describe them. The moduli space of the B-model is the complex
moduli space ${\cal M}_C$ of $M$. Since the holomorphic $(3,0)$-form $\Omega$
is unique up to scale, it defines a line bundle ${\cal L}$ over ${\cal M}_C$
(a sub-bundle of the Hodge bundle) with a natural metric,
\eqn\kahlerpot{ || \Omega ||^2 = i \int_M \Omega \wedge \bar \Omega.}
The metric on ${\cal M}_C$ is given as a curvature of this line bundle, 
\eqn\Kahlermetric{ G_{i \bar j} = \partial_i \bar \partial_{\bar j} K,}
with the K\"ahler potential $K$,
\eqn\complexkahler{K = - {\rm log}|| \Omega ||^2. }

On ${\cal M}_C$, there are particularly useful coordinates called the flat coordinates. 
To define them, let us choose a symplectic basis of homology 3-cycles,
$\{ \alpha_I, \beta^I \}_{I=0,1,\cdots, h^{1,2}}$. Note that 
${\rm dim} H_3 = 2 + 2h^{1,2}$ since $h^{3,0}=1$ for a Calabi-Yau manifold.
Let us consider the period integrals of $\Omega$, 
\eqn\periods{ X^I = \int_{\alpha_I}\Omega , ~~ F_I = \int_{\beta^I}\Omega .}
Since the complex structure of $M$ is determined by the periods $X^I$
over the $\alpha$ cycles, $F_I$'s are functions of 
$X^I$'s. Moreover, they are homogeneous functions of $X$'s 
with weight $1$ since $X^I$'s and $F_I$'s scale in the same 
way under scaling of $\Omega$.
It is known that they can be expressed as derivatives of a single
function $F(X)$ as,
\eqn\prepotential{ F_I(X) = {\partial F(X)\over \partial X^I} , }
where $F(X)$ is a homogeneous function of $X$'s of weight $2$. 
In physics literature, $F(X)$ is called the prepotential. 
Since scaling of $\Omega$ does not affect the complex structure,
we can use ratios of $X$'s
as coordinates of the moduli space ${\cal M}_C$,
\eqn\flatcoordinates{ t^i = {X^i \over X^0},  ~~ i=1 , \cdots, h^{1,2}.}
These $t$'s are called the flat coordinates of ${\cal M}_C$ and play an 
important role in defining the mirror map from $M$ to its mirror 
partner. 
The prepotential $F(X)$ is not globally defined as a
holomorphic section of ${\cal L}^2$ on the moduli space, 
because the periods $X^I, F_I$ undergo monodromy transformations.  
However, its third derivatives with respect to $t$'s,
\eqn\yukawaC{ C_{ijk} = 
{\partial^3  F\over \partial t^i \partial t^j \partial t^k}
= \int_M \Omega {\partial^3 \Omega\over \partial t^i \partial t^j \partial t^k} ,
}
define a global holomorphic section of ${\cal L}^{2}
\otimes {\rm Sym}^{\otimes 3} T^{*1,0} {\cal M}_K$.
In physics literature, $C_{ijk}$ is called the Yukawa couplings.
It is equal to the genus-$0$ topological string amplitude
${\cal A}_{g=0}(3,0)$
evaluated for 3 elements of ${\rm H}^{1}(M, T^{0,1}M)
\simeq {\rm H}^{1,2}(M)$ corresponding to the deformations
$\partial_i, \partial_j, \partial_k$ of complex structure of $M$. 
Since the moduli space of a sphere with three punctures
${\cal M}_{g=0}(3,0)$ is $0$-dimensional, the string amplitude
$\omega_{g=0}(3,0)$ is equal to the conformal field theory
amplitude ${\cal A}_{g=0}(3,0)$. 

\lref\Cecotti{
  S.~Cecotti and C.~Vafa,
  ``Topological antitopological fusion,''
  Nucl.\ Phys.\  B {\bf 367}, 359 (1991).
  %%CITATION = NUPHA,B367,359;%%
}
 
%\HoriIC
\lref\mirrorbook{
K.~Hori, S.~Katz, A.~Klemm, R.~Pandharipande, R.~Thomas, C.~Vafa, R.~Vakil,
and E.~Zaslow, ``Mirror symmetry,'' Clay Mathematics Monographs 1, 
(AMS, Providence, 2003).}

The existence of the line bundle ${\cal L}$ over the moduli space 
and its relation to the moduli space metric are general features of
the topological string theory, not restricted to the B-model \Cecotti . 
The $Q$ cohomology of ${\cal H}$ makes
an analogue of the Hodge bundle over the moduli space and it can
be used to define the metric and compute the genus-$0$ string amplitudes. 
For more details, see Chapter 2 of \BCOV .

In particular, the same structure exists for the A-model, even though 
we must move away from classical geometry in this case. Classically,
the K\"ahler moduli space ${\cal M}_K$ is a cone in $H^{1,1}(M)$ 
where the metric in $M$ is positive definite and the Yukawa coupling
is a constant tensor given by the intersection,
\eqn\classicalyukawa{ C^{({\rm classical})}_{abc} 
= \int_M u_a u_b u_c ,}
of $u_a, u_b, u_c \in {\rm H}^{1,1}(M)$. Quantum mechanically, the metric
and the Yukawa coupling are corrected due to nontrivial holomorphic maps
from $\Sigma_{g=0}$ to $M$. For example, the quantum corrected Yukawa coupling
is given by
\eqn\quantumyukawa{ C_{abc} = 
\int_M u_a u_b u_c  + \sum_n  n_a n_b n_c N_{0,n}{e^{2\pi int} \over 
1 - e^{2\pi i nt}},}
where $N_{0, n}$ is the genus-$0$ Gromov-Witten invariants for
holomorphic maps of degrees $n=(n_1,\cdots,n_{h^{1,1}})$, 
and $t=(t^1, \cdots , t^{h^{1,1}})$
are the flat coordinates of ${\cal M}_K$. In the A-model, $t$'s are
simply the linear coordinates on ${\rm H}^{1,1}(M)$. Once the Yukawa 
coupling is given, one can integrate $C_{abc} = \partial_a \partial_b
\partial_c F$ to compute the prepotantial $F$. The K\"ahler potential
$K$ is then given by    
\eqn\kahlerpot{ K =- {\rm log}\left( 4F - 4 \bar F + \bar t^a \partial_a F
- t^a \bar \partial_{\bar a} \bar F \right).}
Due to the quantum corrections \quantumyukawa , the K\"ahler metric 
$G_{a\bar b} = \partial_a \bar\partial_{\bar b} K$ is not the one for 
the cone of $H^{1,1}(M)$ anymore, but it acquires more elaborate structure. 
\lref\greene{
  P.~S.~Aspinwall, B.~R.~Greene and D.~R.~Morrison,
``Multiple mirror manifolds and topology change in string theory,''
  Phys.\ Lett.\  B {\bf 303}, 249 (1993)
  [arXiv:hep-th/9301043];
  %%CITATION = PHLTA,B303,249;%%
  ``Calabi-Yau moduli space, mirror manifolds and spacetime 
topology change in string theory,''
  Nucl.\ Phys.\  B {\bf 416}, 414 (1994)
  [arXiv:hep-th/9309097].
  %%CITATION = NUPHA,B416,414;%%
}
\lref\wittenphase{  E.~Witten,
``Phases of ${\cal N} = 2$ theories in two dimensions,''
  Nucl.\ Phys.\  B {\bf 403}, 159 (1993)
  [arXiv:hep-th/9301042].
  %%CITATION = NUPHA,B403,159;%%
}
In particular, K\"ahler moduli spaces of topologically distinct Calabi-Yau manifolds
can be combined together into a single smooth moduli space \refs{\wittenphase,
\greene}.

A pair of Calabi-Yau manifolds $(M, \tilde{M})$ is called a mirror pair if 
the A-model for $M$ is equivalent to the B-model for $\tilde M$.
It is an isomorphism of two conformal field theories. Thus, for example,
the Hilbert space ${\cal H}$ of the A-model for $M$ is isomorphism to
the Hilbert space $\tilde {\cal H}$ of the B-model for $\tilde M$, 
including how they are decomposed into irreducible representations 
of the ${\cal N}=2$ superconformal algebra.    
The first example of the mirror pair was found by B.~Greene and M.~R.~Plesser
\lref\GreenePlesser{
  B.~R.~Greene and M.~R.~Plesser,
``Duality in Calabi-Yau moduli space,''
  Nucl.\ Phys.\  B {\bf 338}, 15 (1990).
  %%CITATION = NUPHA,B338,15;%%
} \GreenePlesser .
\lref\Candelas{
  P.~Candelas, X.~C.~De La Ossa, P.~S.~Green and L.~Parkes,
  ``A pair of Calabi-Yau manifolds as an exactly soluble superconformal
  theory,''
  Nucl.\ Phys.\  B {\bf 359}, 21 (1991).
  %%CITATION = NUPHA,B359,21;%%
}
The flat coordinate $t$ that appear in \quantumyukawa\ in the A-model 
are ratios of the periods of $\Omega$ in the B-model and 
the Yukawa coupling \yukawaC\ is also expressed in terms of the
periods. Therefore, 
one can compute the genus-$0$ Gromov-Witten invariants of $M$ by
evaluating the period integrals in $\tilde M$. This procedure was
applied by P.~Candelas, X.~De La Ossa, P.~Green, and L.~Parkes \Candelas\
to the mirror pair of Greene and Plesser. It would be fair to say
that it was their computation that demonstrated the power of 
the mirror symmetry for the first time and sparked interests in
 the mathematical community. 

\newsec{Topological String at Higher Genera}

At genus-$0$, the A-model computes the number of holomorphic maps from 
a sphere to the Calabi-Yau manifold $M$, and the B-model amplitudes are
given by the periods of the holomorphic $(3,0)$-form. The mirror symmetry
relates these two computations. Since the mirror symmetry is an isomorphism 
of two conformal field theories, we expect that the relation 
\eqn\abhigherg{F_g^{\rm A-model}(M)
= F_g^{{\rm B-model}}(\tilde M). }
continues to hold for $g \geq 1$. We note that $F_g$ is a section of
${\cal L}^{2-2g}$.

\subsec{Genus One} 

The genus-$1$ case is special for two reasons \BCOVone. 
Firstly, we do not have to restrict to the
case of $\hat c = 3$ to consider the vacuum amplitude $F_g$. This is because,
for $g=1$, the 
index theorem mentioned in the last section only says that $\eta$
and $\theta$ have the same number of zero modes for $g=1$, and
this does not prevent $\omega_{g=1}(0,0)$ from carrying a top form on 
${\cal M}_{g=1}$ for any $\hat{c}$. 
Secondly, since $\Omega^{0,1}(\Sigma)$ and $\Omega^{1,0}(\Sigma)$ 
are trivial on a genus-$1$ surface, the A-model \fermionsa\ and the B-model
\fermionsb\ have the same set of fields. What happens is that the genus-$1$ vacuum
amplitude is a sum of two contributions, one depends on the K\"ahler moduli
$(t,\bar t)$ and another depends on the complex moduli $(\tau, \bar\tau)$
of $M$,
\eqn\genusonesplit{ F_{1} = F_{1}^{{\rm A-model}}(t,\bar t) +  
F_{1}^{{\rm B-model}}(\tau,\bar \tau).}

It is instructive to consider the case when $M$ is an elliptic curve $T^2$. We use the
standard parameter $\tau$ for the complex moduli of $T^2$. The K\"ahler moduli $t$ can
be chosen so that its imaginary part is the area of $T^2$. In this case, the 
conformal field theory on $\Sigma$ consists of a complex-valued massless free
scalar field (corresponding to the map $X: \Sigma \rightarrow T^2$) and 
a few massless free fermions. These fields are
free since the metric on $T^2$ is flat. This means that their path integrals 
are Gaussian with coefficients given by the Laplacians on $\Sigma$. Thus the conformal
field theory amplitude $\omega_{g=1}$ is given by a sum over {\it harmonic maps} from
$\Sigma$ to $T^2$ weighted by the determinants of the Laplacians on $X$ and on 
$\eta, \theta$. The determinants cancel out due to the $Q$ invariance, and the sum over
harmonic maps and the integral of the resulting $\omega_{g=1}$ over the moduli space
${\cal M}_{g=1}$ can be carried out using the Poisson re-summation method. The result is
\eqn\etafunctions{ F_{1} = - {\rm log}\left(\sqrt{{\rm Im}t}|\eta(t)|^2\right)
- {\rm log}\left(\sqrt{{\rm Im}\tau}|\eta(\tau)|^2\right),}
where $\eta(\tau)$ is the Dedekind eta-function. 
We note that the genus-$1$ amplitude is expressed as a sum of the $t$ dependent term
and the $\tau$ dependent term, as expected in \genusonesplit . On the other hand,
the holomorphic splitting is spoiled by the holomorphic anomalies,
\eqn\torusanomaly{ {\partial^2  F_{1}\over \partial t \partial \bar t}
 = {1 \over 2 ({\rm Im}t)^2}, ~~
{\partial^2  F_{1}\over \partial \tau \partial \bar \tau}
 = {1 \over 2 ({\rm Im}\tau)^2}.}

The A-model part of $F_1$ counts the number of holomorphic maps
in an appropriate sense. To see this, it is useful to write  
\eqn\torussum{ {i \over 2\pi }
{\partial F_{1} \over \partial t}_{|\bar t \rightarrow \infty}
= {1 \over 24} - {1 \over 2} \sum_{M}
e^{2\pi i |{\rm det} M| t},}
where the sum in the right-hand side is over $M \in GL(2,Z)$ 
such that it maps $\tau$ into the fundamental domain of ${\cal M}_{g=1}$.
This can be interpreted as a sum over holomorphic maps
from $\Sigma$ to $T^2$. For a generic complex structure of $\Sigma$, there
is no holomorphic map from $\Sigma$ to $T^2$. Thus, the counting of holomorphic maps
makes sense only when we integrate over ${\cal M}_{g=1}$.
The complex moduli dependent part 
$-{\rm log}(\sqrt{{\rm Im}\tau}|\eta(\tau)|^2)$
has the familiar expression as the logarithm 
of the determinant of the Laplacian on $T^2$. 
It is worthwhile to note
that the torus $T^2$ is self-mirror. Namely, the mirror 
of $T^2$ is another $T^2$ with the K\"ahler moduli and
the complex moduli exchanged. The exchange symmetry
of \etafunctions\ under $t \leftrightarrow \tau$ is 
a consequence of the mirror symmetry in this case. 

For a general Calabi-Yau manifold $M$, the splitting \genusonesplit\ continues to hold.
The A-model amplitude $F_{1}^{{\rm A-model}}$ computes the genus-$1$ Gromov-Witten
invariants,\foot{Note that some of the expressions for $F_{1}$ given in 
\BCOVone\ are missing a factor of $1/2$ since we did not properly take into
account the $Z_2$ automorphism of $T^2$. This has been corrected in \BCOV.} 
\eqn\genusoneGW{
\eqalign{ {i \over 2\pi} 
{\partial F_{1}^{{\rm A-model}} \over \partial t^a}_{|\bar t \rightarrow \infty}  =& 
{(-1)^{\hat c} \over 24} \int_M u_a \wedge c_{\hat c - 1}
- \sum_n  n_a N_{1,n} \sum_m {m e^{2\pi i mnt} \over 1 - e^{2\pi i mnt}}\cr
& - {1 \over 12} \sum_n n_a  N_{0,n} {e^{2\pi in t} \over 1 - e^{2\pi i nt}},}}
where $c_n$ is the $n$-th Chern-Class of $M$, and
$u_a \in H^{1,1}(M)$ corresponds to the deformation
of the K\"ahler moduli $\partial/\partial t^a$. 
Note the contributions from the genus-$0$ 
Gromov-Witten invariants in the second
line in the right-hand side. A mathematical explanation of these contributions is
given by S.~Katz in the appendix to \BCOVone. 

The B-model amplitude $F_{1}^{{\rm B-model}}$ can be expressed in terms of
determinants of del-bar operators on various bundles over $M$.  
Let $V$ be a holomorphic vector bundle over $M$
and $\Delta_V^{(p)} = \bar \partial_V^\dagger \bar \partial_V$ on $\Omega^{0,p}(M)
\otimes  V$. The holomorphic Ray-Sinter torsion for $V$ is defined by 
\eqn\raysinger{ I(V) = \prod_p 
\left( {\rm det}' \Delta_V^{(p)}\right)^{(-1)^p p}.}
It is important to note that $I(V_1)/I(V_2)$ is independent of 
the K\"ahler structure of $M$ \lref\raysinger{D.~B.~Ray and 
I.~M.~Singer, 
``Analytic torsion for complex manifolds,''
Ann. Math. 98, 54 (1973).} \raysinger . 
The genus-$1$ B-model amplitude is given by their combination as
\eqn\raysingerB{ F_{1}^{{\rm B-model}}
 = {1 \over 2} \sum_q (-1)^q q
{\rm log}I(\Omega^{q,0}).}

In \BCOVone , the holomorphic anomaly equation for $F_{1}$ was derived using
the invariance of $\omega_{g=1}$ under $D = d + Q$. The idea is to use the 
invariance to relate the $Q$ trivial operation (such as 
$\partial^2/\partial t \partial \bar t$) on $F_{1}$ to contributions from the
boundary of ${\cal M}_{g=1}$. The result,
for a Calabi-Yau threefold $M$, is given by
\eqn\holgenusone{
\partial_i \bar \partial_{\bar j} F_1 = 
{1 \over 2} C_{imn} \bar C_{\bar j
\bar m \bar n} e^{2K} G^{m\bar m} G^{n\bar n}
-\left( {\chi(M) \over 24} - 1\right) G_{i\bar j},}
where $\chi(M)$ is the Euler characteristic of $M$. In sections 5.6 - 5.8 
of \BCOV, it was shown that this agrees with the application of
the Quillen anomaly formula,
\eqn\quillen{ \partial \bar \partial
{\rm log} I(V) =  
\partial \bar \partial \sum_p
(-1)^p d_p + 2\pi i \int_M {\rm Td}(TM){\rm Ch}(V),}
where $d_p$ is the determinant of the 
inner product in the kernel of $\bar \partial_V$ on $\Omega^{0,p}(M)\wedge V$,
${\rm Td}$ is the Todd class and ${\rm Ch}$ is the Chern class. 
However, the \holgenusone\ does not assume that $F_1$ is given in terms
of the Ray-Singer torsions. In particular, it applies to the A-model also. 

The mirror symmetry \abhigherg\ for the genus-$1$ amplitudes
relates the genus-$1$ Gromov-Witten invariants in \genusoneGW\ to 
the Ray-Singer torsions \raysingerB. In \BCOVone, this was used to compute
$N_{g=1,n}$ explicitly for the quintic threefold. The result was
recently proven mathematically by A.~Zinger 
\lref\zinger{A.~Zinger, ``The reduced
genus-one Gromov-Witten invariants
of Calabi-Yau hypersurfaces,'' [math/0705.2397].} \zinger. 
The result of \BCOVone\ also predicted
 a certain behavior of the torsions \raysingerB\ 
at boundaries of the moduli space of 
the mirror of the quintic threefold. This
prediction was also proven recently
by H.~Fang, Z.~Lu, and 
K.~Yoshikawa \lref\yoshikawa{H.~Fang, Z.~Lu, and K.~Yoshikawa, 
``Analytic torsion for Calabi-Yau threefolds,'' [math.DG/060141].}
\yoshikawa. 

\subsec{Genus greater than One}

The A-model amplitude $F_g^{{\rm A-model}}$ 
for $g > 1$ is related to the genus-$g$ Gromov-Witten
invariants as (sections 5.10 - 5.13 of \BCOV), 
\eqn\Ahigherg{
\eqalign{ F_{g|\bar t \rightarrow \infty}^{{\rm A-model}} 
= &{1 \over 2}\chi(M) \int_{{\cal M}_g}
c^3_{g-1} + \sum_n N_{g,n}e^{2\pi i n t}\cr
&+({\rm contributions~from~}N_{g'<g,n}{\rm ~and~multicoverings}),}} 
where $c_n$ is the $n$-th Chern class of the Hodge bundle over ${\cal M}_g$. 
The formula for the first term in the right-hand side was derived in
\lref\fb{C.~Faber and R.~Pandharipande, ``Hodge integrals and Gromov-Witten
theory,'' math.AG/9810173.} \fb\ as 
\eqn\faber{ \int_{{\cal M}_g} c_{g-1}^3 
= {(-1)^{g-1} \over (2\pi)^{2g-2}} 2 \zeta(2g-2) \chi_g,}
where $\chi_g$ is the Euler characteristics of ${\cal M}_g$ given by
\eqn\eulermg{ \chi_g = {(-1)^{g-1} \over 2g(2g-2)} B_g.}

The B-model amplitude is formally expressed as a sum of Feynman diagrams for
quantization of the Kodaira-Spencer theory (sections 5.1-5.4 of \BCOV).
Although such an expression is a natural generalization of the genus-one 
result relating $F_1$ to the Ray-Singer torsion -- in quantum 
field theory, one-loop amplitudes
are given by determinant of Laplace operators -- not much is known on how 
to carry out the computation beyond $g=1$ using Feynman diagrams.

The holomorphic anomaly equation found in \BCOV\ has turn out to be a useful
tool to compute $F_g$ for higher genera,
\eqn\holanomaly{ \bar\partial_{\bar i} F_g
= {1 \over 2} \bar C_{\bar i \bar j \bar k}e^{2K} G^{j\bar j}
G^{k\bar k}\left( D_j D_k F_{g-1} + \sum_{r=1}^{g-1}
 D_j F_{r} D_k F_{g-r}\right).}
It is the only known method to compute $F_g$ 
systematically for higher genera in the case of
{\it compact} Calabi-Yau 
manifolds. In \BCOV, a diagramatic method was developed
to obtain a general solution to the holomorphic anomaly
equation recursively 
in $g$, to all order in $g$. Recently, the method 
was made more efficient by S.~Yamaguchi and S.-T.~Yau
\lref\Yamaguchi{
  S.~Yamaguchi and S.~T.~Yau,
  ``Topological string partition functions as polynomials,''
  JHEP {\bf 0407}, 047 (2004)
  [arXiv:hep-th/0406078].
  %%CITATION = JHEPA,0407,047;%%
} \Yamaguchi . In \lref\Klemm{
  M.-X.~Huang, A.~Klemm and S.~Quackenbush,
  ``Topological string theory on compact Calabi-Yau: modularity and boundary
  conditions,''
  arXiv:hep-th/0612125.
  %%CITATION = HEP-TH/0612125;%%
} \Klemm , it was shown that, using known behavior of 
$F_g$ at boundaries of the Calabi-Yau moduli space, 
the holomorphic anomaly equations can 
be integrated to obtain $F_g$ up to $g=51$ for the quintic threefold. It is 
hoped that, with a better understanding of boundary data, we can go to even higher genera. 
In fact, the asymptotic behavior of $F_g$ for $g \rightarrow \infty$
for {\it compact} Calabi-Yau manifolds is needed to clarify the relation 
between the
topological string theory and quantum states of black holes to be 
discussed later.  

\newsec{Open/Closed String Duality}

In section 2.2, we defined open string theory for a given choice of
boundary conditions $\{ {\cal E}_i: i = 1, \cdots, k\}$. 
Let us denote the topological string amplitude of genus-$g$
with $n_i$ boundaries of type ${\cal E}_i$ ($i=1,\cdots, k$) 
by $F_{g; n_1,...,n_k}$. It is useful to consider a generating
function with weights $t_1, \cdots , t_k$ as,
\eqn\sumboundary{ F_g(t) = \sum_{n_1,\cdots, n_k} 
F_{g; n_1,...,n_k} t_1^{n_1} \cdots t_k^{n_k}.}
The open/closed string duality is a conjecture that there is a family 
$\tilde M(t)$ of a Calabi-Yau threefold such that the generating 
function $F_g(t)$ defined in the
above is the genus-$g$ vacuum amplitude for {\it closed} topological string
on $\tilde M(t)$. The idea of open/closed string duality 
originally appeared 
in the paper \lref\thooft{G.~'t Hooft, ``A planar diagram theory
for strong interactions,'' Nucl.\ Phys.\ B {\bf 72}, 461 (1974).
%%CITATION = NUPHA,B72,461%%
} by G.~'t Hooft \thooft\ 24 years ago. 
For this reason, $t_i$'s are called 't Hooft couplings. 

Let us present a couple of examples for this. 
The vacuum amplitude of the Chern-Simons gauge theory on $S^3$ 
is given by \lref\wittenjones{
  E.~Witten,
  ``Quantum field theory and the Jones polynomial,''
  Commun.\ Math.\ Phys.\  {\bf 121}, 351 (1989).
  %%CITATION = CMPHA,121,351;%%
}
\wittenjones
\eqn\cs{ Z(S^3)
= {e^{i {\pi \over 8}N(N-1)}
\over (k+N)^{{N\over 2}}}
\sqrt{{k+N \over N}}\prod_{s=1}^{N-1}
\left[ 2 {\rm sin}\left({s\pi \over k+N}\right) \right]^{N-s}.}
Here $k$ is the level of the Chern-Simons theory and the gauge
group is $U(N)$. In \wittenCS , it was observed that $Z(S^3)$ can be expressed
as 
\eqn\cs{ Z(S^3) = {\rm exp}\left( - \sum_{g,n} F_{g, n} \lambda^{2g-2} t^n
\right), }
where the string coupling $\lambda$ and the 't Hooft coupling $t$ are given
in terms of the Chern-Simons variables as
\eqn\thooftexp{ \lambda = {2\pi \over k+N}, ~~ t = i \lambda N,}
and that $F_{g,n}$ is the A-model amplitude of genus $g$
with $n$ boundaries. The target space is the cotangent space of
$S^3$, which is a non-compact Calabi-Yau threefold, and the boundary
condition is set so that open strings end on the base $S^3$ of
$T^* S^3$, or more specifically the map $X :\Sigma_{g,n} \rightarrow 
T^*S^3$ obeys the Neumann condition along $S^3$ and
the Dirichlet condition transverse to $S^3$. The boundary conditions
for the fermions $\eta, \theta$ are determined so that the $Q$
symmetry is preserved. The fact that the base $S^3$ is 
a Lagrangian sub-manifold of 
$T^* S^3$ then guarantees the $Q$ symmetry. 
Physicists refer to this boundary 
condition as ``D branes are wrapping the base $S^3$ of $T^*S^3$.'' 

Subsequently, Gopakumar and Vafa pointed out \lref\Gopakumar{
  R.~Gopakumar and C.~Vafa,
  ``Topological gravity as large $N$ topological gauge theory,''
  Adv.\ Theor.\ Math.\ Phys.\  {\bf 2}, 413 (1998)
  [arXiv:hep-th/9802016].
  %%CITATION = 00203,2,413;%%
} 
that the generating function of $F_{g,n}$ with the 't Hooft coupling
$t$ can be expressed as \Gopakumar\
\eqn\gv{\eqalign{ F_g(t) = &  \sum_n F_{g, n} t^n \cr
   =& \int_{{\cal M}_g} c_{g-1}^3
   - {\chi_g\over (2g-3)!} \sum_{n=1} n^{2g-3} e^{2\pi i n t},}}
where the Hodge integral and the Euler characteristics $\chi_g$ in the
right-hand side are given in the last section in \faber\  and \eulermg .
This is exactly equal to the genus-$g$ A-model amplitude for 
the small resolution of the conifold $u^2 + v^2 + x^2 + y^2 = 0$,
where we can identify $t$ as the K\"ahler moduli for the blown up $P^1$. 
Thus, they conjectured that 
the open string theory on $T^*S^3$ is equivalent to the closed
string theory on the resolved conifold, where the 't Hooft coupling $t$
in the open string side is identified with the amount of blow up on the 
closed string side. 

\lref\OVknots{
  H.~Ooguri and C.~Vafa,
  ``Knot invariants and topological strings,''
  Nucl.\ Phys.\  B {\bf 577}, 419 (2000)
  [arXiv:hep-th/9912123].
  %%CITATION = NUPHA,B577,419;%%
}

Since the Chern-Simons gauge theory can be used to compute invariants of knots 
and links in 
three dimensions \wittenjones ,
it is natural to ask if there are closed 
string duals for such computations also. 
To compute these invariants in the 
topological string theory, one needs to
introduce another D brane ($i.e.$ another Lagrangian sub-manifold of $T^* S^3$)
which intersects with the base $S^3$ along the knot or link in question.
One can take the closed string dual of this set up and find that the same
knot invariants can be computed in the resolved conifold  \OVknots.  
This result predicted new algebraic 
structure of knot and link invariants
and the existence of new integer invariants, 
which were proven in 
\lref\liu{
  K.~Liu and P.~Peng,
  ``Proof of the Labastida-Marino-Ooguri-Vafa conjecture,''
  arXiv:0704.1526 [math.QA].
  %%CITATION = ARXIV:0704.1526;%%
} \liu . 

In this way, knot invariants can be regarded as
topological string amplitudes with
2 types of D branes; one is the base $S^3$ of $T^* S^3$ and another 
is determined by the choice of knot. A generalization of this for 3 types of D branes 
gives the topological vertex, which can be used to compute $F_g$ for any toric 
Calabi-Yau manifold \lref\topvertex{
  M.~Aganagic, A.~Klemm, M.~Marino and C.~Vafa,
  ``The topological vertex,''
  Commun.\ Math.\ Phys.\  {\bf 254}, 425 (2005)
  [arXiv:hep-th/0305132].
  %%CITATION = CMPHA,254,425;%%
} \topvertex. For more detail, see \lref\marinobook{
  M.~Marino,
  ``Chern-Simons theory, matrix models, and topological strings,''
%\href{http://www.slac.stanford.edu/spires/find/hep/www?irn=6927467}{SPIRES entry}
( Oxford, UK: Clarendon ,2005).}
 \marinobook .

So far we have been discussing open topological string associated to the A-model.
According to \wittenCS , an analogue of the Chern-Simons gauge theory in the B-model
is the holomorphic Chern-Simons theory. A particularly simple case is when D branes
wrap $P^1$ in a Calabi-Yau space. Consider the singular space, 
\eqn\singular{ u^2 + v^2 + y^2 + W'(x)^2 = 0 ~~{\rm in } ~~C^4,}
where $W(x)$ is a degree $(k+1)$ polynomial $x^{k+1}+ \cdots$. We can make small resolution
and blow up $k$ $P^1$'s. We can then consider $k$ different boundary conditions,
one for each $P^1$. It was shown in \matrix\ that $F_{g; n_1, \cdots , n_k}$
in this case  is computable by the random matrix model with a potential given by
${\rm tr} W(M)$, 
\eqn\dv{ \int dM  \exp\left[-{1\over \lambda}{\rm tr} W(M)\right]
 = {\rm exp}\left[- \sum_{g, n_1,\cdots ,n_k} F_{g;n_1,\cdots , n_k}\lambda^{2g-2}
t_1^{n_1}\cdots t_k^{n_k}\right],}
where $\int dM$ is an integral over $N \times N$ matrices, $t_i = i\lambda N_i$,
and $N_i$ is the number of eigenvalues of $M$ near the $i$-th critical point of 
$W(x)$. The closed string dual is the topological string on the deformed geometry,
\eqn\deformed{ u^2 + v^2 + y^2 + W'(x)^2 + f(x) = 0, ~~(u,v,x,y) \in C^4,}
for some polynomial $f(x)$ determined by $t_i$'s. 

It is important to note that the open/closed string duality holds for each genus 
separately. As such, it is a property of two-dimensional conformal field theories. 
In \lref\wsderivation{
  H.~Ooguri and C.~Vafa,
``Worldsheet derivation of a large $N$ duality,''
  Nucl.\ Phys.\  B {\bf 641}, 3 (2002)
  [arXiv:hep-th/0205297].
  %%CITATION = NUPHA,B641,3;%%
} \wsderivation, a microscopic derivation of the Gopakumar-Vafa duality was given
by studying phases of the conformal field theory. It is hoped that the AdS/CFT
correspondence \lref\maldacena{
  J.~M.~Maldacena,
  ``The large $N$ limit of superconformal field theories and supergravity,''
  Adv.\ Theor.\ Math.\ Phys.\  {\bf 2}, 231 (1998)
  [Int.\ J.\ Theor.\ Phys.\  {\bf 38}, 1113 (1999)]
  [arXiv:hep-th/9711200].
  %%CITATION = IJTPB,38,1113;%%
} \lref\magoo{
  O.~Aharony, S.~S.~Gubser, J.~M.~Maldacena, H.~Ooguri and Y.~Oz,
  ``Large $N$ field theories, string theory and gravity,''
  Phys.\ Rept.\  {\bf 323}, 183 (2000)
  [arXiv:hep-th/9905111].
  %%CITATION = PRPLC,323,183;%%
} \refs{\maldacena,\magoo} can be derived in a similar way. 

\newsec{Quantum States of Black Holes}

So far we have discussed mainly mathematical aspects of topological string 
theory. Indeed, in the earlier days, the topological string theory was
studied as a toy model of string theory. It is worth pointing out, however, 
that Witten already expressed the vision in his pioneering paper 
\lref\wittensigma{
  E.~Witten,
  ``Topological Sigma Models,''
  Commun.\ Math.\ Phys.\  {\bf 118}, 411 (1988).
  %%CITATION = CMPHA,118,411;%%
} \wittensigma\ that the topological string theory may 
describe a new phase of string theory where general covariance is
unbroken.\foot{In the Minkowski space, general covariance is spontaneously
broken to the Poincare symmetry.} By the early 90's, it was also realized
that the genus-$0$ topological string
amplitude $F_0(M)$ can be used for some computation in physical
superstring compactified on $M$ times the four-dimensional Minkowski 
space. The third derivatives of $F_0$ given by \yukawaC\ and
\quantumyukawa\ are called the Yukawa couplings since they are related to
the couplings of two spinors to one scalar field in the heterotic string
compactified on a Calabi-Yau manifold \lref\distler{
  J.~Distler and B.~R.~Greene,
  ``Some exact results on the superpotential from Calabi-Yau
  compactifications,''
  Nucl.\ Phys.\  B {\bf 309}, 295 (1988).
  %%CITATION = NUPHA,B309,295;%%
} \distler . See also section 6 of 
\wittenmirror.

In \BCOV, the relation to topological string and 
physical string was generalized for all $g\geq 1$. It turned out that 
a certain series of higher derivative terms in the low energy
effective action of type IIA / IIB superstring theory compactified
on $M$ times the Minkowski space is computed by $F_g$'s of the A / B-model.
But, what are these terms in the low energy effective theory good for? 
One answer to this question was found in black holes.

 \lref\OSV{
  H.~Ooguri, A.~Strominger and C.~Vafa,
  ``Black hole attractors and the topological string,''
  Phys.\ Rev.\  D {\bf 70}, 106007 (2004)
  [arXiv:hep-th/0405146].
  %%CITATION = PHRVA,D70,106007;%%
}
%\LopesCardosoWT
\lref\dewit{
  G.~Lopes Cardoso, B.~de Wit and T.~Mohaupt,
  ``Corrections to macroscopic supersymmetric black-hole entropy,''
  Phys.\ Lett.\  B {\bf 451}, 309 (1999)
  [arXiv:hep-th/9812082].
  %%CITATION = PHLTA,B451,309;%%
}

Black holes are classical solutions to the Einstein equation coupled to matter
fields. There are singularities, but they are covered by event horizons.
Anything that goes inside of the horizons will not be able to come back out. 
As classical solutions, black holes are characterized by a very small set of
parameters such as mass, angular momentum, and electric and magnetic charges. 
In contrast, a quantum black hole is expected to carry a large number of
quantum states, roughly the exponential of the area of the horizon measured
in the unit of the Planck area $\sim 10^{-70}{\rm  m}^2$. For example, 
a typical astrophysical black hole of mass $\sim 10^{31}{\rm kg}$ would have about 
$e^{10^{78}}$ states. This expectation is based on the following. 

For classical solutions to the Einstein equation with an energy-momentum tensor
satisfying the positive energy condition, S.~Hawking proved the black hole area
theorem stating that the total area of black hole event horizons cannot 
decrease in {\it any} classical process 
\lref\areatheorem{S.~W.~Hawking, ``Gravitational radiation from
colliding black holes,'' Phys.\ Rev.\ Lett.\ {\bf 26},
1344 (1971)} \areatheorem.  Based on an analogy between this theorem 
and the second
law of thermodynamics, which states that the thermodynamic
entropy can only increase in time, 
J.~Bekenstein proposed that a black hole carries 
an entropy $S$ proportional to the area $A$ of its horizon \lref\Bekenstein{
  J.~D.~Bekenstein,
  ``Black holes and entropy,''
  Phys.\ Rev.\  D {\bf 7}, 2333 (1973).
  %%CITATION = PHRVA,D7,2333;%%
} \Bekenstein . The proportionality constant was undetermined at that time,
however. Within a year, S.~Hawking  \lref\Hawking{
  S.~W.~Hawking,
  ``Black hole explosions,''
  Nature {\bf 248}, 30 (1974).
  %%CITATION = NATUA,248,30;%%
} \Hawking\
used the semi-classical quantization of spacetime geometry and matter
near the horizon and showed that a black hole of mass $M$ emits the 
black body radiation with temperature $T$ consistent with the 
standard thermodynamic formula
\eqn\temp{ {1\over T } = {\partial S \over \partial M},}
if the entropy $S$ is given by 
\eqn\bhentropy{ S = {1\over 4} {A \over \ell_{{\rm P}}^2},}
where $\ell_{{\rm P}} \sim 1.6 \times 10^{-35} {\rm m}$ is the Planck length.
Hawking's computation gave an evidence for Bekenstein's conjecture
and also fixed the proportionality constant to be $(4 \ell_{{\rm P}}^2)^{-1}$. 
Since the area $A$ of the horizon scales as $M^2$, this entropy formula
predicts that a black hole of mass $M$ should have quantum 
states as many as $\sim e^{\left({M\over M_{{\rm P}}}\right)^2}$, 
where the Planck mass $M_{{\rm P}} \sim 
10^{-8}{\rm kg}$. The number $e^{10^{78}}$ quoted in the previous paragraph
is obtained by estimating $e^S$ for astrophysical black holes. 
The proposed entropy formula \bhentropy\ can be compared to the definition
of the thermodynamic entropy by R.~Clausius in the 19th century. 
The Bekenstein-Hawking formula \bhentropy\ was proposed based on 
macroscopic properties of black holes, 
and it calls for a microscopic explanation just as Boltzmann 
recognized that the thermodynamic 
entropy is the logarithm of the number of states and explained much
of classical thermodynamics based on the microscopic point of view. 
Doing the same for the black hole entropy was posed as a challenge to any theory
that claims to unify quantum mechanics and general relativity consistently. 

String theory has partially met this challenge. There is a class of black holes
in superstring theory called BPS black holes. These solutions are
invariant under non-trivial subalgebras of supersymmetry, and their Hawking 
temperatures are zero. These black holes are therefore stable, 
and their quantum states can be reliably counted. Thanks to the D brane
construction by Polchinski
\lref\Polchinski{
  J.~Polchinski,
  ``Dirichlet-branes and Ramond-Ramond charges,''
  Phys.\ Rev.\ Lett.\  {\bf 75}, 4724 (1995)
  [arXiv:hep-th/9510017].
  %%CITATION = PRLTA,75,4724;%%
} \Polchinski , some of them can be realized in terms of D branes in superstring
theory. This description is useful since low energy quantum states
of these black holes can be described using supersymmetric gauge 
theories on D branes. 
A.~Strominger and C.~Vafa used this technique to count the number of quantum states of 
BPS black holes and found a perfect agreement with the prediction of Bekenstein
and Hawking in the limit when the masses of black holes are asymptotically 
large \lref\strominger{
  A.~Strominger and C.~Vafa,
  ``Microscopic Origin of the Bekenstein-Hawking Entropy,''
  Phys.\ Lett.\  B {\bf 379}, 99 (1996)
  [arXiv:hep-th/9601029].
  %%CITATION = PHLTA,B379,99;%%
} \strominger .
One has to take this limit since Hawking's computation in \Hawking\ is reliable
only in this limit. As the black hole becomes small, the curvature near the horizon
becomes strong and the semi-classical approximation breaks down. 

To deepen our understanding of the black hole microstates, we need to know
what happens for small black holes. In fact, the counting of microstates typically
becomes simpler for small black holes since the gauge theory on D branes
becomes weakly coupled. On the other hand, the gravity side of the story 
becomes
more complicated.
 The strong curvature at the horizon means that we must take into
account corrections to Einstein's theory due to stringy effects and quantum gravity
effects. Remarkably, topological string amplitudes $F_g$ compute exactly such 
corrections for small BPS black holes. Inspired by the earlier work \dewit , the 
following conjecture was formulated in \OSV . 

For definiteness, let us consider type IIB superstring theory compactified on 
a Calabi-Yau manifold. For type IIA, we take
the mirror of the story below. Let us choose the symplectic basis of homology 3-cycles,
$\{ \alpha_I, \beta^I \}_{I=0,1,\cdots, h^{1,2}}$ on $M$. In type IIB string, there are 
D3 branes -- it means that we can consider boundary conditions such that open strings
end on a $(3+1)$ dimensional submanifold in the four-dimensional Minkowski 
space times $M$. Submanifols are chosen to be special Lagrangian submanifolds of 
$M$, as required by supersymmetry 
\lref\Becker{
  K.~Becker, M.~Becker and A.~Strominger,
  ``Five-branes, membranes and nonperturbative string theory,''
  Nucl.\ Phys.\  B {\bf 456}, 130 (1995)
  [arXiv:hep-th/9507158].
  %%CITATION = NUPHA,B456,130;%%
}
\lref\OOZ{
  H.~Ooguri, Y.~Oz and Z.~Yin,
  ``D-branes on Calabi-Yau spaces and their mirrors,''
  Nucl.\ Phys.\  B {\bf 477}, 407 (1996)
  [arXiv:hep-th/9606112].
  %%CITATION = NUPHA,B477,407;%%
}
\refs{\Becker,\OOZ}, 
times a
time-like direction in the Minkowski space. Seen in four dimensions, it is
a particle moving in the time-like direction. We can consider a bound state of 
several D branes. Take a bound state of D3 branes whose total homology class
is given by $\sum_I p^I \alpha_I + q_I \beta^I$ for some integers $(p^I, q_I)$. 
Seen in four dimensions, the mass of the resulting particle is proportional to the
total volume of the special Lagrangian submanifolds. Clearly, the mass of the particle
becomes large as $(p^I, q_I)$ increases. As the mass becomes large, the particle
starts influencing the spacetime around it. In the limit of large mass,
the Schwarzschild radius becomes longer than the Compton wave-length, 
and the geometry
is described by BPS black hole solutions of the type discussed in the previous
paragraph. What we want to compute is the number of quantum states of D3 branes
for arbitrary values of $(p^I, q_I)$. Let us call the number $\sigma(p,q)$.

To state the conjecture on $\sigma(p,q)$, we need to set up some notation for topological
string amplitudes $F_g$. In the B-model, $F_g$ is a section of ${\cal L}^{2-2g}$
over the moduli space of complex structure of $M$. The complex structure can
be parametrized by the periods $X^I = \int_{\alpha_I} \Omega$ of the holomorphic
$(3,0)$-form. Since ${\cal L}$ is the sub-bundle of the Hodge bundle associated
to scaling of $\Omega$, we can regard $F_g$ as a homogeneous function of $X^I$'s
of weight $(2-2g)$. 

The conjecture of \OSV\ is that
\eqn\osv{ \sum_{q_I} \sigma(p,q) e^{- \pi\phi^I q_I}
 = \Big| {\rm exp}\big( - F(X) \big) \Big|^2, }
where $F(X)$ is the sum of the topological string amplitudes to all genera, 
\eqn\allgenera{
F(X) = \sum_{g=0}^\infty F_g(X), }
and the periods $X^I$ are fixed as 
\eqn\attractor{ X^I = p^I + i \phi^I. }
Since $F_g(X)$ is a homogeneous function of $X$ of weight $(2-2g)$, 
we should think of $F(X)=\sum_{g=0}^\infty F_g(X)$ 
in \allgenera\ as an asymptotic expansion for  
large $X$. If we neglect the higher genus terms $F_{g \geq 1}(X)$ in
\allgenera\ and simply 
use the leading term $F_0(X)$ for $F(X)$ in \osv, 
we reproduce the Bekenstein-Hawking formula \bhentropy\
by $S={\rm log}\sigma(p,q)$ evaluated for large $(p,q)$. 
Thus, the conjecture is that the subleading 
terms \allgenera\ -- terms with $g \geq 1$ -- represent
quantum gravity corrections to the entropy formula. 

The gauge theory computation of $\sigma(p,q)$ takes a familiar form in
type IIA superstring theory. In this case, instead of D3 branes wrapping 3
cycles in $M$, we should consider D0, D2, D4 and D6 branes on holomorphic
cycles in $M$. In particular, D0 branes are points and D6 branes
wraps the entire Calabi-Yau manifold.  
By the homological mirror symmetry
of Kontsevich, $p^I$ count the numbers of D4 and D6 branes, and $q_I$
count the numbers of D0 and D2 branes. If there are no D6 branes,
the gauge theory on D4 branes is the ${\cal N}=4$ supersymmetric gauge theory 
{\it topologically twisted} 
and the number of quantum states is computed by the Witten index 
as in the work by Vafa and Witten \lref\VafaWitten{
  C.~Vafa and E.~Witten,
  ``A Strong coupling test of S duality,''
  Nucl.\ Phys.\  B {\bf 431}, 3 (1994)
  [arXiv:hep-th/9408074].
  %%CITATION = NUPHA,B431,3;%%
} \VafaWitten . It is the Euler characteristics of the instanton moduli space.
The numbers of D0 and D2 branes are the second and the first Chern classes 
of the gauge bundles on D4 branes. The left-hand side of  
\osv\ then becomes the generating function of the 
Euler characteristics of the instanton moduli space on D4 branes
of the type studied by H.~Nakajima \lref\nakajima{H.~Nakajima, ``Instantons on 
ALE spaces, quiver varieties, and Kac-Moody
algebras,'' Duke Math. {\bf 76}, 365 (1994).} \nakajima .
With D6 branes, one has to think about a gauge theory in six dimensions.

\lref\donaldson{
S.~K.~Donaldson and R.~P.~Thomas, ``Gauge theory in higher dimensions,''
in the Geometric Universe (Oxford Univ. Press, 1998).}

Therefore, in mathematical terms, the conjectured formula \osv\ 
relates the Gromov-Witten
invariants of $M$ to the Euler characteristics of 
the moduli space of instantons on four manifolds embedded in
$M$, in the case of the A-model. 
For some {\it non-compact} Calabi-Yau manifolds, both sides
of the formula can be evaluated explicitly, resulting in 
non-trivial checks of the conjecture 
\lref\vafatorus{
  C.~Vafa,
  ``Two dimensional Yang-Mills, black holes and topological strings,''
  arXiv:hep-th/0406058.
  %%CITATION = HEP-TH/0406058;%%
}
%\AganagicJS
\lref\aganagic{
  M.~Aganagic, H.~Ooguri, N.~Saulina and C.~Vafa,
  ``Black holes, $q$-deformed 2d Yang-Mills, and non-perturbative topological
  strings,''
  Nucl.\ Phys.\  B {\bf 715}, 304 (2005)
  [arXiv:hep-th/0411280].
  %%CITATION = NUPHA,B715,304;%%
} \refs{\vafatorus,\aganagic}. Though the formula \osv\ has not been
proven even to the standard of physicists, there have been several
promising attempts and they have deepened our understanding of topological
string, black holes microstates, and the AdS/CFT correspondence
\lref\mooreproof{
  F.~Denef and G.~W.~Moore,
  ``Split states, entropy enigmas, holes and halos,''
  arXiv:hep-th/0702146.
  %%CITATION = HEP-TH/0702146;%%
}
\lref\deboerproof{
  J.~de Boer, M.~C.~N.~Cheng, R.~Dijkgraaf, J.~Manschot and E.~Verlinde,
  ``A farey tail for attractor black holes,''
  JHEP {\bf 0611}, 024 (2006)
  [arXiv:hep-th/0608059].
  %%CITATION = JHEPA,0611,024;%%
}
\lref\dijkgraafproof{
  R.~Dijkgraaf, C.~Vafa and E.~Verlinde,
  ``M-theory and a topological string duality,''
  arXiv:hep-th/0602087.
  %%CITATION = HEP-TH/0602087;%%
}
%\BeasleyUS
\lref\stromingerproof{
  C.~Beasley, D.~Gaiotto, M.~Guica, L.~Huang, A.~Strominger and X.~Yin,
  ``Why $Z_{{\rm BH}} 
= |Z_{{\rm top}}|^2$,''
  arXiv:hep-th/0608021.
  %%CITATION = HEP-TH/0608021;%%
} 
\lref\gaiottproof{
  D.~Gaiotto, A.~Strominger and X.~Yin,
  ``From AdS$_3$/CFT$_2$ to black holes / topological strings,''
  JHEP {\bf 0709}, 050 (2007)
  [arXiv:hep-th/0602046].
  %%CITATION = JHEPA,0709,050;%%
}
\lref\dabholkarproof{
  A.~Dabholkar, F.~Denef, G.~W.~Moore and B.~Pioline,
  ``Precision counting of small black holes,''
  JHEP {\bf 0510}, 096 (2005)
  [arXiv:hep-th/0507014].
  %%CITATION = JHEPA,0510,096;%%
}
\lref\senproof{
  R.~K.~Gupta and A.~Sen,
  ``AdS$_3$/CFT$_2$ to AdS$_2$/CFT$_1$,''
  arXiv:0806.0053 [hep-th].
  %%CITATION = ARXIV:0806.0053;%%
}
\refs{\dabholkarproof,\gaiottproof, \dijkgraafproof, 
\stromingerproof, \deboerproof, \mooreproof, \senproof}.
To understand \osv\ better, it appears useful to develop a method to  
estimate of $F_g$ for large $g$ for {\it compact} Calabi-Yau manifolds. 

\newsec{Concluding Remarks} 

There are many aspects of the topological string theory I could not cover
in this set of lectures. In particular, the topological vertex \topvertex\
to compute $F_g$ for toric Calabi-Yau manifolds is only briefly mentioned
in section 5. Among other important discoveries I could not mention
in these lectures are
the relations of the topological string theory
to dimer models and crystal melting
\lref\crystal{
  A.~Okounkov, N.~Reshetikhin and C.~Vafa,
  ``Quantum Calabi-Yau and classical crystals,''
  arXiv:hep-th/0309208.
  %%CITATION = HEP-TH/0309208;%%
} \lref\crystaltwo{
  A.~Iqbal, N.~Nekrasov, A.~Okounkov and C.~Vafa,
  ``Quantum foam and topological strings,''
  JHEP {\bf 0804}, 011 (2008)
  [arXiv:hep-th/0312022].
  %%CITATION = JHEPA,0804,011;%%
} \crystal, to the Donaldson-Thomas theory 
\lref\moop{D.~Maulik, A.~Oblomkov, A.~Okounkov, and R.~Pandharipande,
``Gromov-Witten/Donaldson-Thomas correspondence for toric 3-folds,''
arXiv:math.AG/0809.3976.}
\refs{\donaldson, 
\crystaltwo,\moop},
to the Seiberg-Witten theory 
\lref\nekrasov{
  N.~Nekrasov and A.~Okounkov,
  ``Seiberg-Witten theory and random partitions,''
  arXiv:hep-th/0306238.
  %%CITATION = HEP-TH/0306238;%%
} \nekrasov , and to BPS counting in M theory 
\lref\Gopatwo{
  R.~Gopakumar and C.~Vafa,
  ``M-theory and topological strings. II,''
  arXiv:hep-th/9812127.
  %%CITATION = HEP-TH/9812127;%%
}
%\GopakumarII
\lref\Gopaone{
  R.~Gopakumar and C.~Vafa,
  ``M-theory and topological strings. I,''
  arXiv:hep-th/9809187.
  %%CITATION = HEP-TH/9809187;%%
} \refs{\Gopaone, \Gopatwo}. More recently,
the non-commutative version of the Donaldson-Thomas
theory has been formulated in
\lref\sz{
  B.~Szendr\"oi,
  ``Non-commutative Donaldson-Thomas theory and the conifold,''
  Geom.\ Topol.\  {\bf 12}, 1171 (2008)
  [arXiv:0705.3419 [math.AG]].
  %%CITATION = 00256,12,1171;%%
}
\lref\mr{
  S.~Mozgovoy and M.~Reineke,
  ``On the noncommutative Donaldson-Thomas invariants 
arising from brane tilings,''  arXiv:0809.0117 [math.AG].
}
\refs{\sz,\mr} and its relation to the crystal melting model
and the black hole microstate counting have been 
pointed out in 
\lref\oy{
  H.~Ooguri and M.~Yamazaki,
  ``Crystal Melting and Toric Calabi-Yau Manifolds,''
  arXiv:0811.2801 [hep-th].
  %%CITATION = ARXIV:0811.2801;%%
} \oy\ generalizing the result of \crystal .

Though the topological string theory has rich mathematical 
structure and a broad range of applications to problems in 
physics, one should be reminded that it is still a
simplified version of superstring theory. There is a vast
{\it terra incognita} in superstring theory, and we need 
new mathematics to guide us through this territory. 
In this connection, it is also humbling to note that,
after 80 years since the conception of 
 quantum field theory by Heisenberg and Pauli 
\lref\Heisenberg{
  W.~Heisenberg and W.~Pauli,
  ``Zur Quantendynamik der Wellenfelder,''
  Z.\ Phys.\  {\bf 56}, 1 (1929).
  %%CITATION = ZEPYA,56,1;%%
} \Heisenberg , there has been no proof that it exists as
a consistent mathematical theory, except for some cases in lower dimensions. 
In fact, the existence proof of quantum version of the Yang-Mills theory
in four dimensions with a demonstration of its confinement property
is posed as one of the Seven Millennium Problems by the Clay
Mathematics Institute. 

After the discovery of D brane construction by 
Polchinski and the AdS/CFT correspondence by Maldacena, 
many quantum field theories have been found to be dual ($i.e.$, equivalent)
 to string theory in various geometries. This duality between quantum
field theory and string theory has been used to evaluate non-perturbative
effects in quantum field theories that are not accessible in any
other method. At the same time, many fundamental issues in string theory
have been translated into quantum field theory questions and this has shed
important lights on mysteries of quantum gravity such
as the black hole information paradox of Hawking. Because of this
development, research in quantum field theories
has become increasingly intertwined with string theory research.
The progress in the past 10 years suggests that an ultimate
solution to the Yang-Mills Problem could be built on
a foundation where quantum field theories and string theory
are extended and transformed into a single mathematical framework.

\bigskip
\medskip
\centerline{\bf Acknowledgments} 

I would like to thank the organizing committee of the Takagi Lectures for
inviting me to present these lectures and for their hospitality in Kyoto. 
I would like to thank  
collaborators on my works described in these lectures, especially 
M.~Bershadsky, S.~Cecotti, A.~Strominger, and C.Vafa.

This work is supported in part by 
DOE grant DE-FG03-92-ER40701,
by a Grant-in-Aid for Scientific 
Research (C) 20540256 from the Japan Society 
for the Promotion of Science,
by the World Premier International Research Center 
Initiative of MEXT of Japan, and by the Kavli Foundation.

\listrefs
\end